\documentclass[a4paper,final,10pt]{article}
\usepackage{subfigure}
\usepackage[colorinlistoftodos]{todonotes}
\usepackage{xcolor}
\usepackage{amsmath}
\usepackage{amssymb}
\usepackage{bm}
\usepackage{ulem}

\newtheorem{theorem}{Theorem}[section]

\newtheorem{lemma}[theorem]{Lemma}
\newtheorem{proposition}[theorem]{Proposition}
\newtheorem{definition}[theorem]{Definition}

\newtheorem{remark}[theorem]{Remark}
\newtheorem{conjecture}[theorem]{Conjecture}

\newenvironment{proof}{\textit{Proof:}}{$\hfill\square$\newline}

\newcommand\rev[1]{\textcolor{black}{#1}}
\newcommand\redout{\bgroup\markoverwith{\textcolor{red}{\rule[.5ex]{2pt}{0.8pt}}}\ULon}

\title{Saddle-node canard cycles in \\ planar piecewise linear differential systems}

\author{V. Carmona\footnotemark[1]
\and S. Fern\'andez-Garc\'ia\footnotemark[2]
\and A. E. Teruel\footnotemark[3]}

\begin{document}

\maketitle
\renewcommand{\thefootnote}{\fnsymbol{footnote}}
\footnotetext{\textit{MSC2010 subject classification}: Primary: 34C05, 34C23, 34C25, 34E15, 34E17; Secondary: 37G15, 37G25.}
\footnotetext[1]
{Dpto. Matem\'atica Aplicada II \& IMUS, University of Seville, Escuela Superior de Ingenieros, Avenida de los Descubrimientos s/n, 41092 Sevilla, Spain}
\footnotetext[2]
{Dpto. EDAN \& IMUS, University of Seville, Facultad de Matem\'aticas, C/ Tarfia, s/n., 41012 Sevilla, Spain. Corresponding author: soledad@us.es}
\footnotetext[3]
{Departament de Matem\`atiques i Inform\`atica \& IAC3, Universitat de les Illes Balears, Palma de Mallorca, Spain}

\begin{abstract}
By applying a singular perturbation approach, canard limit cycles exhibited by a general family of singularly perturbed planar piecewise linear (PWL) differential systems are analyzed. The performed study involves both hyperbolic and non-hyperbolic canard limit cycles appearing after both a supercritical and a subcritical Hopf bifurcation.

The obtained results are completely comparable with those obtained for smooth vector fields. In some sense, the manuscript can be understood as an extension towards the PWL framework of the results obtained for smooth systems by Krupa and Szmolyan \cite{KS01}. In addition, some novel slow-fast behaviors are  obtained. In particular, in the supercritical case, and under suitable conditions, it is proved that the limit cycles are organized along a curve exhibiting two folds. Each of these folds corresponds to a saddle-node bifurcation of canard limit cycles, one involving headless canard cycles, whereas the other involving canard cycles with head. This configuration allows the coexistence of three canard limit cycles. 
\end{abstract}

\section{Introduction}

Planar slow-fast systems are differential systems involving two variables which evolve with different velocities. The ratio of the velocities of both variables is the so-called singular parameter. From a qualitative point of view, 
slow-fast systems exhibit neither remarkable properties, nor special bifurcations.
Phase portraits of slow-fast systems are topologically equivalent to phase portraits of differential systems without different time scales, and canard limit cycles are topologically equivalent to usual limit cycles. The main difference between slow-fast systems and those which are not, is the rate of variation of interesting quantitative information versus the variation of some parameters, specially when the singular parameter is small. This is the case of the phenomenon called canard explosion, where the amplitude of a limit cycle, born at a supercritical Hopf bifurcation, increases very rapidly while a system parameter varies in an exponentially small range, see \cite{dumortier,eckhaus,KS01}.

The canard explosion phenomenon was discovered and analyzed by Benoit et al. in 1981 \cite{benoit} in the Van der Pol oscillator and explains the fast transition, by varying a parameter, from a small amplitude limit cycle to the relaxation oscillation appearing in this system.  Relaxation oscillations are oscillatory behavior characterized by long periods of quasi-static behavior interspersed  with short periods of rapid transition. Since this behavior is usual in real-life phenomena, see \cite {KS01} and references therein, models based in slow-fast differential equations are ubiquitous in many applications, such as chemical and biological ones, and in particular in neuroscience \cite{FN1,izhikevich,FN2}. 

The main tools for the analysis of the slow-fast dynamics are provided by Geometric Singular Perturbation Theory, and rely on the ability of reconstructing the global dynamics by splitting and then joining, in a suitable way, the fast and slow behaviors. Under hyperbolicity conditions, Fenichel Theorem  describes the existence of invariant slow manifolds close to compact parts of the fast nullcline and also describes the stability properties of these slow manifolds \cite{fenichel}. When the fast nullcline folds, the flow of big parts of the phase space evolves in a narrow neighborhood of the stable branch of the slow manifold, and close to the fold. This fact has dynamical consequences since the existence of limit cycles can be deduced from the behavior of the slow manifold around the fold. In some sense, global dynamical behaviors, like limit cycles, can be deduced from local aspects as the local transition next to the fold. 

A usual technique to analyze this local transition is the blow up of the fold \cite{dumortier,KS00,KS01}. Some other authors have analyzed this transition by reproducing it in different contexts, which result more friendly for analysis. Piecewise smooth systems, both continuous and discontinuous, have received the main attention of the authors, and many progress have been made. For instance, in \cite{R16,glend} the authors analyze the canard explosion in a continuous piecewise smooth context. Nevertheless, the simplest scenario where this phenomenon can be reproduced is the piecewise linear (PWL) framework. 
Even when some dynamical aspects of the slow-fast behavior had been observed in PWL systems, see \cite{DFK16,DGPPRT16} and references therein, it has taken some time to understand the way of reproducing the slow-fast dynamics properly, see \cite{chapter,DFHPT13,DGPPRT16}. In \cite{FDKT} the authors reproduce part of the canard explosion phenomenon in the PWL context, in particular the one involving hyperbolic headless canards.        

In this work, we consider an extension of the system analyzed in \cite{FDKT}, that allows for the existence of both canards with and without head and both, hyperbolic and  non-hyperbolic canard cycles. In particular, the system is able to reproduce saddle-node bifurcations of canard limit cycles. 
The obtained results are completely comparable with those obtained, for smooth vector fields, by Krupa and Szmolyan in \cite{KS01}.
Moreover, we find new scenarios that, as far as we are concerned, have not been previously reported in the smooth framework. In particular, we find situations where two saddle-node bifurcations of canard cycles take place, one of headless canards and another one of canards with head. In such a case, we show the coexistence of three canard limit cycles.

The article is organized as follows. In Section \ref{sec:background}, we provide a brief overview of canard explosion and saddle-node canard cycles in the smooth case. After that, in Section \ref{sec:pwlsystem}, we introduce the class of systems we aim to study. In Section \ref{sec_mainresults}, we present the Main Results of the article.  Section \ref{sec_proofs} is devoted to the proofs of the Main Results.
Finally, Section \ref{sec_conclusions} is devoted to conclusions and possible extensions of the present work. The technical issues of the proofs have been left to the Appendix.

\section{Background on canard cycles: canard explosion}\label{sec:background}
In this section, based on \cite{KS00, KS01}, we briefly review the basic ingredients of the  canard oscillatory behavior appearing in the smooth framework. Typically, canard solutions take place in planar differential systems of the form
\begin{equation}\label{smooth1}
\left\{\begin{array}{l}
\varepsilon\dot x=f(x,y,a,\varepsilon), \\
\noalign{\medskip}
~~\dot y= g(x,y,a,\varepsilon),
\end{array}
\right.
\end{equation}
where $f,g\in \mathcal{C}^r, r\geq 3, a\in\mathbb{R}, 0<\varepsilon\ll 1$ and the dot denotes the derivative with respect to the temporal variable $\tau.$ Since the velocity of the solutions of system \eqref{smooth1} is very different depending on the regions of the phase plane they are crossing through (far from the $x$-nullcline the velocity is fast, and close to the $x$-nullcline it is slow),  system \eqref{smooth1} is often called {\it slow-fast system}. 

After the rescaling in time $t=\tau/\varepsilon$, system \eqref{smooth1} writes as
\begin{equation}\label{smooth2}
\left\{\begin{array}{l}
x'=f(x,y,a,\varepsilon), \\
\noalign{\medskip}
y'= \varepsilon g(x,y,a,\varepsilon),
\end{array}
\right.
\end{equation}
where the prime denotes the derivative with respect to the fast time $t.$
Systems (\ref{smooth1}) and (\ref{smooth2}) are equivalent through the identity when $\varepsilon>0$, but they have not the same limit for $\varepsilon=0$. In fact, the limit of system (\ref{smooth1}), called {\it slow subsystem}, is a semi-explicit differential algebraic equation (DAE), where the relation between the variables is given by
\[
 S=\{(x,y):f(x,y,a,0)=0\}. 
\]
Assuming that $f_y(x,y,a,0)\neq 0$ it follows that $S$ is the graph of a differentiable function $y=\varphi_a(x)$, and the DAE reduces to the differential equation 
\begin{eqnarray}\label{eq:reduced}
 f_x(x,\varphi_a(x),a,0)\dot{x}=-f_y(x,\varphi_a(x),a,0)g(x,\varphi_a(x),a,0),
\end{eqnarray}
which is called the {\it reduced equation}. On the other hand, the limit for $\varepsilon=0$ of system (\ref{smooth2}), called {\it fast subsystem}, is a differential equation having $S$ as the locus of every equilibrium point. From here $S$ is called the {\it critical manifold}.

From Fenichel's Theorems, we obtain an approximation of the overall slow dynamics through the slow subsystem, and an approximation of the overall fast dynamics through the fast subsystem. In brief, Fenichel's theory asserts that compact subsets $S_0 \subset S$ formed by normally hyperbolic equilibrium points of the critical manifold, persist as locally invariant slow manifolds $S_{\varepsilon}$ for $\varepsilon>0,$ which can be extended by the flow. Moreover, $S_{\varepsilon}$ is an attracting or repelling manifold depending on the same character, under the flow defined by the fast subsystem, of $S_0$. Furthermore, the flow over $S_{\varepsilon}$ is a regular perturbation of the flow defined by the slow subsystem, or equivalently the reduced equation. 

Typically, the breakdown of the normal hyperbolicity takes place at the points $(x_0,y_0)\in S$ at which the manifold folds, i.e., $f_x(x_0,y_0,a,0)=0$ and $f_{xx}(x_0,y_0,a,0)\neq0$.  It is possible to assume, without loss of generality, that the fold point is at the origin when $a= 0$, in that case we avoid the subindex at function $\varphi_a$, and therefore $\varphi(x)=\eta x^2+O(x^3)$. Hence, in a neighborhood of the origin two different branches of the critical manifold coexist, the attracting one, $ S^a_0=\{(x,y): f_x(x,y,0,0)<0\}$, and the repelling one, $S^r_0=\{(x,y):f_x(x,y,0,0)>0\}$. Over each of these two branches the reduced equation can be written as the EDO
\[
 \dot{x} = -\frac{f_y(x,\varphi(x),0,0)}{f_x(x,\varphi(x),0,0)}g(x,\varphi(x),0,0)=
 \frac {g(x,\varphi(x),0,0)}{\varphi'(x)},
\]
and, assuming that $g(0,0,0,0)\neq 0$ , the flow defined by it has opposite orientation over each of these branches. In this case the fold point is called {\it jump point}.

A special situation occurs when the fold point at the origin satisfies that $g(0,0,0,0)=0$ together with the non-degeneracy condition $g_x(0,0,0,0)\neq 0$. The reduced equation can be then regularized, defining a solution of the desingularized system which passes from one branch, $S_0^a$, to the other, $S_0^r$, through the fold point. In this case the fold point is called a {\it canard point}. 

\begin{figure}[ht]
\begin{center}
\includegraphics[width=6.5cm]{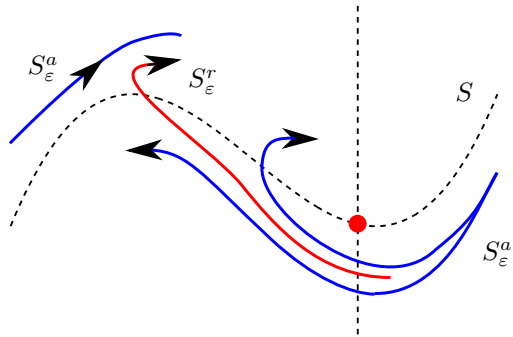}
\begin{picture}(0,0)
\put(-10,30){$S_{\varepsilon}^a$}
\put(-180,100){$S_{\varepsilon}^a$}
\put(-20,90){$S$}
\put(-120,95){$S_{\varepsilon}^r$}
\end{picture}
\end{center}
\caption{Representation of the critical manifold $S$ in the neighborhood of a canard point. Attracting branch, $S_{\varepsilon}^a$, and repelling branch, $S_{\varepsilon}^r$, of the slow manifold obtained after singular perturbation are also represented. The canard cycles with and without head are obtained provided $S_{\varepsilon}^a$ flows along one side or the other of the repelling manifold $S_{\varepsilon}^r.$
}\label{fig:intro}
\end{figure}

After perturbation, i.e., for $\varepsilon>0$, Fenichel slow manifolds $S^a_{\varepsilon}$ and $S^r_{\varepsilon}$ behave in a different way near a jump fold point and near a canard point. Around a jump point, an attracting Fenichel slow manifold $S^a_\varepsilon$ may follow closely the attracting branch $S_0^a,$ pass in the vicinity of the fold
point, and continue following approximately the fast dynamics, giving rise to the possibility of relaxation oscillations. However, around a canard point,  Fenichel slow manifold $S^a_\varepsilon$ may follow closely the attracting branch $S_0^a,$ pass in the vicinity of the fold point and then, surprisingly, continue following closely the repelling branch $S_0^r$, see Figure \ref{fig:intro}. From this behavior it can be concluded the existence of solutions of system \eqref{smooth1} with $0<\varepsilon\ll 1$ containing canard segments. 

Under the existence of another attracting Fenichel slow manifold, see Figure \ref{fig:intro}, when $S^a_\varepsilon$ flows along one side or the other of the  $S^r_\varepsilon$,  they can exist \textit{canards without head} or {\it canards with head}, respectively. Since slow manifolds are exponentially close to one another, the presence of exponentially small terms in the expansions in power series of $\varepsilon$ of the slow manifolds implies that their respective position can change upon an exponentially small parameter variation. This phenomenon is known as the {\it canard explosion}. Moreover, the transition from canards without head to canards with head occurs typically when $S^a_\varepsilon$ connects to $S^r_\varepsilon$. This connection takes place along a curve ${a}_c(\varepsilon)$ in the parameter plane $(\varepsilon,a)$ and the associated canard solution is said to be a \textit{maximal canard}. 

Canard cycles develop along a branch born at a Hopf bifurcation, at $a=a_H$, and the canard explosion takes place at a value which is at a distance of $O(\varepsilon)$ from the $a_H$. This means that very close to the bifurcation point $a_H$, before the explosion, the cycles have the characteristics of typical Hopf cycles. This Hopf bifurcation arises only for $\varepsilon>0$ and is usually known as a singular Hopf bifurcation \cite{braaksma, guckenheimer}.

The existence of saddle-node bifurcation of canard cycles in the smooth framework has been analyzed in \cite{KS01}. There, the authors consider two different cases, depending whether the Hopf bifurcation where the cycle is born is supercritical or subcritical. Thus, after proving the existence of the maximal canard, 
they distinguish two different scenarios: 
\begin{itemize}
\item \textit{Supercritical case}: In Theorem 3.3, authors state the existence of a family of periodic orbits. These periodic orbits can be stable Hopf-type limit cycles, canard limit cycles or relaxation oscillations. To analyze the stability of the canard limit cycles, they use the \textit{way in-way out function} $R(s)$, which is the limit of the integral of the divergence along the slow manifolds when $\varepsilon\to 0$. In Theorem 3.4, assuming that this function is negative, the authors state that the canard limit cycles of the family are stable.
\item \textit{Subcritical case}: In Theorem 3.5, authors state the existence of other family of periodic orbits. The orbits of that family can be unstable Hopf-type limit cycles, canard limit cycles or relaxation oscillations. Again, to analyze the stability of canard cycles, they use the  \textit{way in-way out function} $R(s)$. In Theorem 3.6, assuming that this function has exactly one simple zero at $s=s_{lp,0}$, the authors state that there exists a function $s_{lp}(\sqrt{\varepsilon})$ having limiting point at $s_{lp,0}$ when $\varepsilon\to 0$, such that canard limit cycles are unstable for $s<s_{lp}(\sqrt{\varepsilon})$ and stable for $s>s_{lp}(\sqrt{\varepsilon}).$
 \end{itemize}

Instead of the way in-way out function $R(s)$ used in \cite{KS01}, other authors use a  similar concept known as slow divergence integral, see \cite{D11}.

\section{Statement of the piecewise linear system.}\label{sec:pwlsystem}
In this section we introduce the family of PWL differential systems we are going to work with, together with some basic elements of their dynamics. We also define some functions and quantities which are needed for stating the main results in the next section. 

Let us consider the following family of planar differential systems depending on the four dimensional parameter $\bm{\eta}=(a,k,m,\varepsilon)$,
\begin{equation}\label{systempal}
\left\{\begin{array}{l}
x'=y-f(x,a,k,m,\varepsilon),\\
y'=\varepsilon(a-x),
\end{array}\right.
\end{equation}
where the prime denotes the derivative with respect to the time $t$, $(x,y)^T\in\mathbb{R}^2$, $0<\varepsilon\ll 1,$ and the $x$-nullcline is defined by the graph of the continuous PWL function with four segments given by
\begin{equation}\label{vectorfield4z}
f(x,a,k,m,\varepsilon)=
\left\{
  \begin{array}{ll}
    x+1-k(\sqrt{\varepsilon}-1)-m(\sqrt{\varepsilon}+a),&\quad \mbox{if }x<-1\\
    -k(x+\sqrt{\varepsilon})-m(\sqrt{\varepsilon}+a),	&\quad \mbox{if }-1<x\leq-\sqrt\varepsilon,\\
    m(x-a),						&\quad \mbox{if } |x|\leq\sqrt\varepsilon,\\
    x-\sqrt\varepsilon+m(\sqrt\varepsilon-a), 		&\quad \mbox{if } x>\sqrt\varepsilon,
\end{array}\right.
\end{equation}
with $k>0$ and $|m|<2\sqrt{\varepsilon}$. 

The phase space is splitted into four regions: the lateral half-planes $LL=\{(x,y): x\leq -1\}$ and $R=\{(x,y): x\geq \sqrt{\varepsilon}\}$, and the central bands $L=\{(x,y): -1 \leq x\leq -\sqrt{\varepsilon}\}$ and $C=\{(x,y): |x|\leq\sqrt{\varepsilon}\}$. Restricted to any of these regions, the vector field is linear and it can be expressed in a matrix way as ${F}_i(\mathbf{x}) =A_i\mathbf{x} +\mathbf{b}_i$ with $i\in\{LL,L,C,R\}$ and
\[
 A_{LL}=\left(
    \begin{array}{cc}
      -1 & 1 \\
      -\varepsilon & 0
    \end{array}
 \right),\ 
 A_L=\left(
    \begin{array}{cc}
      k & 1 \\
      -\varepsilon & 0
    \end{array}
 \right),\ 
 A_C=\left(
    \begin{array}{cc}
      -m & 1 \\
      -\varepsilon & 0
    \end{array}
 \right),\ 
 A_R=\left(
    \begin{array}{cc}
      -1 & 1 \\
      -\varepsilon & 0
    \end{array}
 \right),
\]
\[
 \mathbf{b}_{LL}=\left(
    \begin{array}{c}
      (k+m)\sqrt{\varepsilon}+ma-(k+1) \\
      \varepsilon a
    \end{array}
 \right),\ 
 \mathbf{b}_L=\left(
    \begin{array}{c}
      (k+m)\sqrt{\varepsilon}+ma \\
      \varepsilon a
    \end{array}
 \right),\ 
  \mathbf{b}_C=\left(
    \begin{array}{c}
      ma \\
      \varepsilon a
    \end{array}
 \right),
\]
\[
  \mathbf{b}_R=\left(
    \begin{array}{c}
      \sqrt{\varepsilon}(1-m)+ma \\
      \varepsilon a
    \end{array}
 \right).
\]

The local behavior of the flow of system \eqref{systempal}-\eqref{vectorfield4z} at any of the regions $i\in\{LL,L,C,R\}$ is determined by the trace $t_i$, the determinant $d_i=\varepsilon$, the discriminant $\Delta_i=t_i^2-4\varepsilon$, the slow eigenvalue $\lambda_i^s$, the fast eigenvalue $\lambda_i^q$, the slow eigenvector $\mathbf{v}_i^s=(\lambda_i^s,-\varepsilon)^T$ and the fast eigenvector $\mathbf{v}_i^q=(\lambda_i^q,-\varepsilon)^T$ of the matrix $A_i$, and by the location of the point $\mathbf{e}_i=-A_i^{-1}\mathbf{b}_i$.  We summarize all this information in Table \ref{tab_1} and Table \ref{tab_2}.
\begin{table}[h]
\begin{center}
\begin{tabular}[h]{|c|c|c|}\hline 
  & $LL$ &  $R$ \\ \hline
$t_i $ & $-1$ &  $-1$ \\ \hline  
$\Delta_i$ & $1-4\varepsilon$ &  $1-4\varepsilon$ \\ \hline
$\lambda_i^s$ & $\frac {-1+\sqrt{1-4\varepsilon}}{2}=-\varepsilon- \varepsilon^2 + \ldots$ & $\frac {-1+\sqrt{1-4\varepsilon}}{2}=-{\varepsilon}- {\varepsilon^2} + \ldots$ \\ \hline
$\lambda_i^q$ & $-1-\lambda_i^s$ & $-1-\lambda_i^s$\\ \hline
$\mathbf{e}_i$ & $\left(\begin{array}{c} a \\ 1+k-\sqrt{\varepsilon}(m+k)-a(m-1) \end{array}\right)$ &  $\left(\begin{array}{c} a \\ (m-1)(\sqrt{\varepsilon}-a)\end{array} \right)$\\ \hline
\end{tabular}
\end{center}
\caption{Significant quantities for the dynamics of system \eqref{systempal}-\eqref{vectorfield4z} in the lateral half-planes $LL$ and $R$. Point $\mathbf{e}_i=-A_i^{-1}\mathbf{b}_i$ with $i\in\{LL,R\}$ is an equilibrium point only when $\mathbf{e}_{LL}\in LL$ or  $\mathbf{e}_{R}\in R$.}\label{tab_1}
\end{table}

We remark that $\mathbf{e}_i$ are equilibrium points only when they are located in the region where the system \eqref{systempal}-\eqref{vectorfield4z} behaves as the linear system $\mathbf{F}_i(\mathbf{x})=A_i\mathbf{x}+\mathbf{b}_i$. Otherwise, these points are called {\it virtual equilibrium points}, and they also organise the dynamic behaviour in that region even when they are not equilibrium points.

\begin{table}[bh]
\begin{center}
\begin{tabular}{|c|c|c|}\hline 
  &  $L$ & $C$  \\ \hline
$t_i $ &  $k$ & $-m$ \\ \hline  
$\Delta_i$ & $k^2-4\varepsilon$ & $m^2-4\varepsilon$  \\ \hline
$\lambda_i^s$ & $\frac {k-\sqrt{k^2-4\varepsilon}}{2}=\frac{\varepsilon}{k}+\frac {\varepsilon^2}{k^3} + \ldots$ & $-\frac{m}{2} \pm \frac {\sqrt{4\varepsilon-m^2}}{2} i$ \\ \hline
$\lambda_i^q$ & $k-\lambda_i^s$ &  \\ \hline
$\mathbf{e}_i$ & $\left(\begin{array}{c} a \\ -(m+k)(\sqrt{\varepsilon}+a)\end{array}\right)$ & $\left(\begin{array}{c} a \\ 0 \end{array}\right)$ \\ \hline
\end{tabular}
\end{center}
\caption{Significant quantities for the dynamics of system \eqref{systempal}-\eqref{vectorfield4z} in the central bands $L$ and $C$. Point $\mathbf{e}_i=-A_i^{-1}\mathbf{b}_i$ is an equilibrium point only when $\mathbf{e}_i$ is contained in its own region, that is, $\mathbf{e}_{L}\in L$ or  $\mathbf{e}_{C}\in C$.}\label{tab_2}
\end{table}

\begin{figure}[t]
 \includegraphics{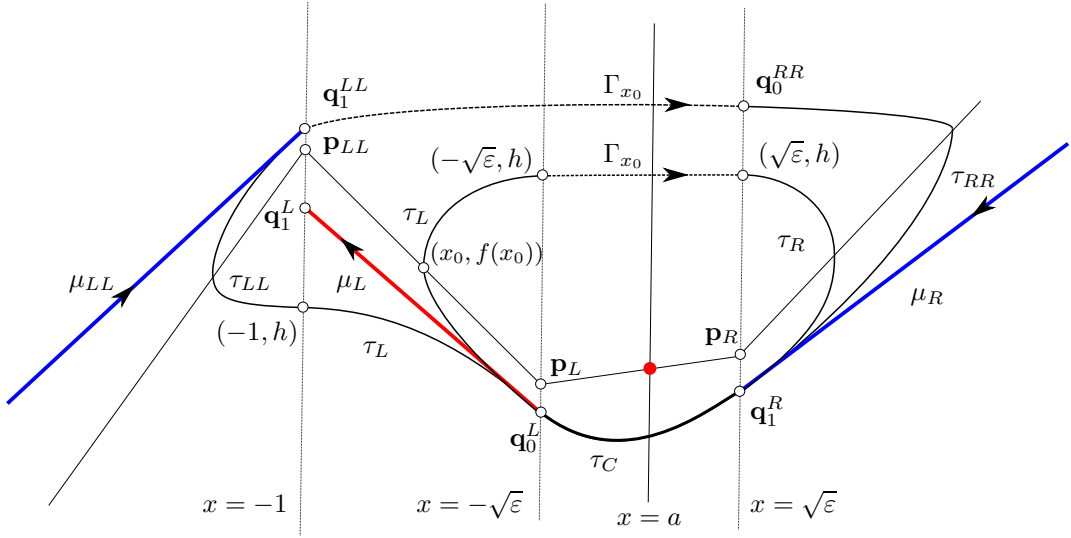}
 \put(-210,35){$\mathbf{q}_0^L$}
 \put(-195,63){$\mathbf{p}_L$}
 \put(-240,140){$(-\sqrt{\varepsilon},h)$}
 \put(-250,120){$\tau_L$}
 \put(-240,105){{\small $(x_0,f(x_0))$}}
 \put(-265,70){$\tau_L$}
 \put(-315,95){$\tau_{LL}$}
 \put(-120,45){$\mathbf{q}_1^R$}
 \put(-137,75){$\mathbf{p}_{R}$}
 \put(-117,142){$(\sqrt{\varepsilon},h)$}
 \put(-110,110){$\tau_R$}
 \put(-175,168){$\Gamma_{x_0}$}
 \put(-175,142){$\Gamma_{x_0}$}
 \put(-180,27){$\tau_C$}
 \put(-170,5){$x=a$}
 \put(-320,75){$(-1,h)$}
 \put(-45,134){$\tau_{RR}$}
 \put(-280,165){$\mathbf{q}_1^{LL}$}
 \put(-280,147){$\mathbf{p}_{LL}$}
 \put(-302,119){$\mathbf{q}_1^{L}$}
 \put(-118,170){$\mathbf{q}_0^{RR}$}
 \put(-375,95){$\mu_{LL}$}
 \put(-275,95){$\mu_{L}$}
 \put(-60,90){$\mu_{R}$}
 \put(-325,10){$x=-1$}
 \put(-245,10){$x=-\sqrt{\varepsilon}$}
 \put(-120,10){$x=\sqrt{\varepsilon }$}
 \caption{Representation of the dynamical objects of system \eqref{systempal}-\eqref{vectorfield4z}. Linearity regions $LL,L,C$ and $R$ and switching lines $x=-1, x=-\sqrt{\varepsilon}$ and $x=\sqrt{\varepsilon}$. The fast nullcline given by the graph of the function $y=f(x,a,k,m,\varepsilon)$  and the intersection points with the switching lines $\mathbf{p}_{LL},\mathbf{p}_L$ and $\mathbf{p}_R$. The slow nullcline $x=a$ and the equilibrium point at the intersection with the fast nullcline.  Slow manifold $S_{\varepsilon}$ with the attracting branch $S_{\varepsilon}^a=\mu_{LL} \cup \mu_{R}$, the repelling branch $S_{\varepsilon}^r=\mu_{L}$, and the intersection points with the switching lines $\mathbf{q}_1^{LL}, \mathbf{q}_1^R$ and $\mathbf{q}_0^L, \mathbf{q}_1^L$, respectively. Limit cycles $\Gamma_{x_0}$ are also represented both cycles with head $x_0<-1$ and cycles without head $x_0\in(-1,-\sqrt{\varepsilon})$.}\label{dib:geommetry}
\end{figure}

From Lemma 4 in \cite{PTV16}, the slow manifold $S_{\varepsilon}$ of system \eqref{systempal}-\eqref{vectorfield4z}, with $0<\varepsilon\ll 1$, is locally formed by segments, each of them contained in a region $i\in \{LL,L,R\}$ and defined by the slow eigenvector $\mathbf{v}_i^s=(\lambda_i^s,-\varepsilon)^T$ associated to the slow eigenvalue $\lambda_i^s$. Hence
\begin{equation}\label{def:Se}
 S_{\varepsilon}=\left\{
    \begin{array}{ll}
      \mu_{LL}=\mathbf{e}_{LL}+r\mathbf{v}_{LL}^s & r\in \left[-\frac {1+a}{\lambda_{LL}^s},+\infty \right),\\
      \noalign{\smallskip}
      \mu_{L}=\mathbf{e}_{L}-r\mathbf{v}_{L}^s & r\in\left[\frac{\sqrt{\varepsilon}+a}{\lambda_L^s},\frac{1+a}{\lambda_L^s}\right],\\  
      \noalign{\smallskip}
      \mu_{R}=\mathbf{e}_{R}-r\mathbf{v}_{R}^s & r\in\left[\frac {a-\sqrt{\varepsilon}}{\lambda_R^s},+\infty\right).\\      
    \end{array}
 \right.
\end{equation}
Since $|m|<2\sqrt{\varepsilon}$, the matrix $A_C$ has complex eigenvalues with modulus equal to $\sqrt{\varepsilon}$, see Table~\ref{tab_2}. Therefore, there is not neither a fast eigenvector nor a slow  eigenvector. Thus, the slow manifold in the central region is not a segment but it is formed by two pieces of curve, the repelling one connecting with $\mu_L$ and the attracting one connecting with $\mu_R$. In Figure \ref{dib:geommetry} both pieces have been represented by a unique curve. Finally,  since there is not a real splitting between fast and slow behavior at the level of the eigenvalues, it follows that neither the repelling effect nor the attracting effect of the slow manifold in the central region increases drastically as $\varepsilon$ tends to zero. This phenomena was introduced in \cite{DFHPT13} and limit cycles flowing close to this kind of slow manifold are called quasi-canards. For this reason, we do not consider this part of the slow manifold as a part neither of the attracting branch, nor the repelling branch. 

We conclude that $S_{\varepsilon}^a=\mu_{LL} \cup \mu_R$ and $S_{\varepsilon}^r=\mu_L$ are the attracting branch and the repelling branch, respectively, of the slow manifold $S_{\varepsilon}$. Moreover, the attracting branch $S_{\varepsilon}^a$ intersects with the switching lines $x=-1$ and $x=\sqrt{\varepsilon}$ at the points
\begin{equation}\label{q1R}
  \mathbf{q}_1^{LL}= \begin{pmatrix}
                      -1\\-\lambda_{LL}^s(1+a)-k(\sqrt{\varepsilon}-1)-m(\sqrt{\varepsilon}+a)
                     \end{pmatrix},\quad
  \mathbf{q}_1^R=\begin{pmatrix}
		      \sqrt{\varepsilon}\\ (m+\lambda_R^s)(\sqrt{\varepsilon}-a)
		 \end{pmatrix},
\end{equation}
respectively, see Figure \ref{dib:geommetry}, whereas the repelling branch $S_{\varepsilon}^r$ intersects the switching lines $x=-1$ and $x=-\sqrt{\varepsilon}$ at the points 
\begin{equation*}\label{q0L}
  \mathbf{q}_0^{L}=\begin{pmatrix}
			-\sqrt{\varepsilon}\\-(m+\lambda_L^s)(\sqrt{\varepsilon}+a)  
                   \end{pmatrix},\quad
  \mathbf{q}_1^L=\begin{pmatrix}
		 -1\\-(m+k)(\sqrt{\varepsilon}+a)+(1+a)\lambda_L^q
		 \end{pmatrix},
\end{equation*}
respectively. We also highlight the intersection points of the $x$-nullcline with the switching lines $x=-1,x=-\sqrt{\varepsilon}$ and $x=\sqrt{\varepsilon}$,  
\begin{equation}\label{pLpLL}
 \mathbf{p}_{LL}=\begin{pmatrix}
                  -1\\k(1-\sqrt{\varepsilon})-m(\sqrt{\varepsilon}+a)
                 \end{pmatrix},\quad
 \mathbf{p}_L=\begin{pmatrix}
               -\sqrt{\varepsilon}\\-m(\sqrt{\varepsilon}+a)
              \end{pmatrix},\quad 
 \mathbf{p}_R=\begin{pmatrix}
               \sqrt{\varepsilon}\\m(\sqrt{\varepsilon}-a)
              \end{pmatrix},
\end{equation}
respectively. Note that the flow at these points is tangent to the switching line. 

Regarding the invariant sets under the flow of system \eqref{systempal}-\eqref{vectorfield4z}, in the following result we show that all of them are included into the rhomboid $\mathcal{R}$ which is limited by the straight lines $\mathbf{e}_{LL} +r\mathbf{v}_{LL}^s$, $\mathbf{e}_{LL} +r\mathbf{e}_{1}$,
$\mathbf{e}_R +r\mathbf{v}_R^s$, and $\mathbf{e}_R +r\mathbf{e}_1$, where $r\in \mathbb{R}$ and $\mathbf{e}_1^T=(1,0)$. The result follows straightforward by analysing the orientation of the flow over the boundary of $\mathcal{R}$.

\begin{lemma}\label{lem:inv_reg}
\begin{itemize}
\item [a)] The rhomboid $\mathcal{R}$, previously defined, is positively invariant under the flow defined by the system \eqref{systempal}-\eqref{vectorfield4z}.
\item [b)] The rhomboid $\mathcal{R}$, previously defined, contains the equilibrium point and every periodic orbit of the system   \eqref{systempal}-\eqref{vectorfield4z}.
\end{itemize}
\end{lemma}

Every limit cycle $\Gamma$ of system \eqref{systempal}-\eqref{vectorfield4z} intersects the $x$-nullcline $(x,f(x,a,k,m,\varepsilon))$ at exactly one point $(x_{\Gamma},f(x_{\Gamma},a,k,m,\varepsilon))$ with $x_{\Gamma}<a$. 
From now on, we call \textit{width of the limit cycle $\Gamma$}, to the first coordinate of this intersection point, that is $x_{\Gamma}$.  

One special limit cycle, assuming that it exists, is the one having width $x=-1$. Such a limit cycle is tangent to the switching line $\{x=-1\}$ at the point $\mathbf{p}_{LL}$, and therefore, it is the separation cycle between the limit cycles intersecting the lateral region $LL$  and those that do not intersect it. In a similar way, the limit cycle having width $x=-\sqrt{\varepsilon}$ is tangent at $\mathbf{p}_L$ to the switching line $\{x=-\sqrt{\varepsilon}\}$ and it is the separation cycle between the limit cycles intersecting the region $L$ and those that do not intersect it. 

When $\varepsilon$ is small enough,
the limit cycles with width $-1<x<-\sqrt{\varepsilon}$ will be referred to as headless canard limit cycles whereas limit cycles with width $x<-1$ will be referred to as canard limit cycles with head. Therefore, the limit cycle with width $x=-1$ will be referred as {\it the transitory canard}, see  \cite{DDR14}, and it is the boundary between headless canard cycles and canard cycles with head.

In addition, every limit cycle with head, $\Gamma_{x_0}$  where $x_0<-1$, intersects the separation line $\{x=-1\}$ at two points. Let $(-1,h)$ be the one below the point $\mathbf{p}_{LL}$, see Figure \ref{dib:geommetry}. Moreover, every headless limit cycle, $\Gamma_{x_0}$ with $x_0\in(-1,-\sqrt{\varepsilon})$, intersects the separation line $\{x=-\sqrt{\varepsilon}\}$ also at two points. Let $(-\sqrt{\varepsilon},h)$ be the one over the point $\mathbf{p}_{L}$. We referred to $h$ as the {\it height} of the limit cycle  $\Gamma_{x_0}$ both in the case with head and in the case without head.

Therefore, any limit cycle can be labeled by the two different quantities that we have denoted by the width and by the height. Let $\Phi$ be the piecewise function which maps the width of a limit cycle into its height, i.e.  
\begin{equation*}
\Phi(x)=\left\{
    \begin{array}{ll}
        \Phi_{3z}(x) & \text{if } x\in[-1,-\sqrt{\varepsilon}),  \\
        \Phi_{4z}(x) & \text{if } x<-1,
    \end{array}
\right.
\end{equation*}
where $\Phi_{3z}$ is defined by the flow of the linear system in the region $L$, and $\Phi_{4z}$ is defined by the flow, in backward time, of the linear system in the region $LL$. Therefore, through $h=\Phi(x)$ we can pass from the width $x$ of a limit cycle $\Gamma$ to its height $h$.
Typically, the height $h$ is more convenient for computational purposes, whereas the width $x$ is more convenient for stating the results.

To analyze the stability of the canard limit cycles in the PWL framework, it is not possible to use the same approach that it is used in the smooth context, since the singularity of the reduced flow can not be removed through the desingularization process, \cite{KS01}. In fact, the reduced equation \eqref{eq:reduced} associated to the system \eqref{systempal}-\eqref{vectorfield4z} writes as
\[
 \dot{x}=
 \left\{
    \begin{array}{ll}
     \frac{x}{k} & x < 0,\\
     -x      & x> 0,
    \end{array}
 \right.
\]
whose desingularized flow is discontinuous at $x=0$. Therefore, the analysis of the stability of limit cycles through techniques based on the singular flow, such as the way in-way out function, \cite{KS00,KS01}, or the slow divergence integral, see \cite{DDR14} and references there in, can not be successfully applied in this context. Nevertheless, this analysis can be performed directly when $\varepsilon>0$ by explicitly computing the integral of the divergence as the sum of the products of the traces of the linear systems and the time of flight in each region of linearity \cite{freire1}.

\section{Statement of the Main Results}\label{sec_mainresults}

In this section we present the main results in the paper. These results concern to the existence of a one parameter family of canard limit cycles in the PWL system \eqref{systempal}-\eqref{vectorfield4z}, and to the description about how this family organizes along a curve in the plane $(x,a)$, where $x$ is the width of the canard limit cycle and $a$ is the parameter value. The results also provide information about the stability of the limit cycles, paying special attention to semi-stable ones. In order to be fluid in the exposition, we left the proofs and their technical issues  for next sections.

In the first result we assure that, the starting point of the curve organizing the family of limit cycles exhibited by system \eqref{systempal}-\eqref{vectorfield4z} takes place at a Hopf-like bifurcation \cite{freire1}. At this bifurcation a limit cycle appears after the change of stability of the singular point, just like in the Hopf bifurcation. The difference between both kind of bifurcations is the relation between the amplitude of the limit cycle and the bifurcation value, this relation is linear in the Hopf-like bifurcation and a square root in the Hopf bifurcation. 

The proof of the following result is a straightforward conclusion of Theorem 5.1 and Theorem 5.2 in \cite{S19}, see also \cite{freire1} and \cite{S18}.

\begin{theorem} \label{th:hopf}
System \eqref{systempal}-\eqref{vectorfield4z} has a unique singular point $\mathbf{e}=(a,f(a))$ which converges to the fold of the critical manifold at the origin as $(\varepsilon,a)$ tends to zero. Moreover, the function, 
\[
a_H(\sqrt{\varepsilon})=
\left\{
  \begin{array}{ll}
   \sqrt{\varepsilon}& m=-\sqrt{\varepsilon},\\
   -\sqrt{\varepsilon} & m=\sqrt{\varepsilon},
  \end{array}
\right.
\] 
satisfies that the equilibrium is stable for $a>a_H(\sqrt{\varepsilon})$ and looses stability through a Hopf-like bifurcation as $a$ passes through $a_H(\sqrt{\varepsilon}).$  In particular, if $m=-\sqrt{\varepsilon}$, a stable limit cycle appears when $a<a_H(\sqrt{\varepsilon})$ in a supercritical bifurcation, and if $m=\sqrt{\varepsilon}$, a unstable limit cycle appears when $a>a_H(\sqrt{\varepsilon})$ in a  subcritical bifurcation. In both cases, the size of the limit cycle depends linearly on the distance $|a_{H}(\sqrt{\varepsilon})-a|$. 
\end{theorem}

Next theorem is devoted to the existence of a trajectory connecting the attracting branch and the repelling branch of the slow manifold.  This connection is usually referred to as the { maximal canard} trajectory.

\begin{theorem}\label{th:connection}
Set $m=\pm\sqrt{\varepsilon}.$  There exist a value $\varepsilon_0>0$ and a function $a=\tilde{a}(k,\varepsilon; m),$ analytic as a function of $(k,\sqrt{\varepsilon})$, defined in the open set $U=(0,+\infty)\times (0,\varepsilon_0)$ and such that, for $(k,\varepsilon)\in U$,
a solution of system \eqref{systempal}-\eqref{vectorfield4z} starting in the attracting branch of the slow manifold, $\mu_R$, connects to the repelling branch of the slow manifold, $\mu_L$, if and only if $a=\tilde{a}(k,\varepsilon;m)$. In such case, the time of flight of the transition is $\tau_C(k,\varepsilon; m)>0$. First terms of the expansions of $\tilde{a}(k,\varepsilon; m)$ and $\tau_C(k,\varepsilon; m)$   are given as follows,

\begin{equation}\label{aepmenos}
 \tilde{a}(k,\varepsilon;m)=
 \left\{\begin{array}{l}
 \dfrac{ {{e}^{\frac{\pi }{\sqrt{3}}}}-1}{{{e}^{\frac{\pi }{\sqrt{3}}}}+1}\sqrt{\varepsilon}
 -\dfrac{{{e}^{\frac{\pi }{\sqrt{3}}}}}{{{\left( {{e}^{\frac{\pi }{\sqrt{3}}}}+1\right) }^{2}}}
   \left(\dfrac{1-k^2}{{k}^{2}}\right) {{\varepsilon}^{{3}/{2}}}
  +O(\varepsilon^2),\quad \mbox{if }m=-\sqrt{\varepsilon}, \\
-\dfrac {e^{\frac{\pi}{\sqrt{3}}}-1} {e^{\frac{\pi}{\sqrt{3}}}+1}\sqrt{\varepsilon}
 -\dfrac {e^{\frac{\pi}{\sqrt{3}}}} {\left(e^{\frac{\pi}{\sqrt{3}}}+1\right)^2} 
  \left( \dfrac {1-k^2}{k^2} \right) \varepsilon^{3/2}+O(\varepsilon^2),\quad \mbox{if }m=\sqrt{\varepsilon}, 
\end{array}\right.
\end{equation}
and
\begin{equation}\label{tauepmenos}
 \tau_C(k,\varepsilon; m)=
 \left\{\begin{array}{l}
 \dfrac{2 \pi }{\sqrt{3}}\dfrac 1{\sqrt{\varepsilon}}
  -\dfrac{1+k}{k}
  -\dfrac{1-k^2}{2 {{k}^{2}}} \sqrt{\varepsilon}
  +O(\varepsilon),\quad \mbox{if }m=-\sqrt{\varepsilon}, \\ \noalign{\smallskip}
\dfrac{2\pi}{\sqrt{3}} \dfrac1{\sqrt{\varepsilon}}
 -\dfrac {1+k}{k} +  \dfrac {1-k^2}{2k^2} \sqrt{\varepsilon} +O(\varepsilon)  ,\quad \mbox{if }m=\sqrt{\varepsilon}.
\end{array}\right.
\end{equation}
\end{theorem}

The existence of the maximal canard trajectory, together with the divergence of the flow in a neighborhood of the slow manifold, provide the arguments we use in Section \ref{sec_proofs} to proof the following result about the existence of canard cycles of any suitable width. To state the result in a proper way we introduce the following values
\begin{equation}\label{def:xs}
 x_r= -(1+k)+k\sqrt{\varepsilon}-\lambda_L^s (\sqrt{\varepsilon}+a), \quad  
 x_s= -\sqrt{\varepsilon}-\lambda_L^s (\sqrt{\varepsilon}+a).
\end{equation}
This values correspond with the end points of the interval such that limit cycles having width contained in $(x_r,x_s)$ are canard limit cycles. In fact, limit cycles having width $x<x_r$ are relaxation oscillations whereas limit cycles having width $x>x_s$ are still under the effect of the Hopf-like bifurcation.  

\begin{theorem}\label{th:exist_Gamma}
Fix $\varepsilon_0$ sufficiently small and set $m=\pm\sqrt{\varepsilon}$. 
There exists a function $a=\hat{a}(k,\varepsilon,x_0; m),$ $C^\infty$ function of $(k,\sqrt{\varepsilon},x_0)$, defined in the open set $U=(0,+\infty)\times (0,\varepsilon_0)\times(x_r,x_s)$,  fulfilling  
\[
\begin{array}{ll}
| \hat{a}(k,\varepsilon,x_0; m)-\tilde{a}(k,\varepsilon;m)| \approx |x_0| e^{-\frac {x_0}{\varepsilon^{3/2}}} & x_0\in[-1,x_s),\\
|\hat{a}(k,\varepsilon,x_0; m)-\tilde{a}(k,\varepsilon;m)| \approx |x_0-x_r| e^{-\frac {x_0-x_r}{\varepsilon}} & x_0\in(x_r,-1),
\end{array}
\]
with $\tilde{a}(k,\varepsilon;m)$ the function defined in Theorem \ref{th:connection}, and such that, for $(k,\varepsilon,x_0)\in U$ and $a=\hat{a}(k,\varepsilon,x_0; m)$ system \eqref{systempal}-\eqref{vectorfield4z} possesses a canard limit cycle, $ \Gamma_{x_0}$, passing through $(x_0,f(x_0))$. The canard limit cycle is headless if $x_0\in(-1,{x}_s)$ and with head if $x_0\in(x_r,-1).$ 
\end{theorem}

Previous result describe the canard explosion taking place in the PWL framework. There, it can be observed that the slope of the explosion is different before and after the maximal canard. 

\begin{figure}[ht]
\begin{center}
 \includegraphics[width=4.5cm]{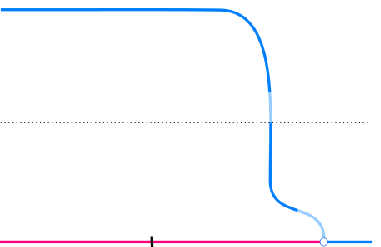}\quad \includegraphics[width=4.5cm]{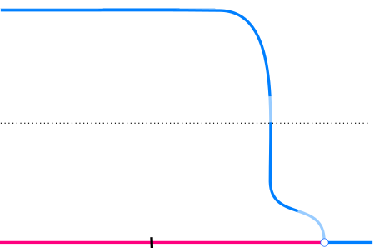}\quad \includegraphics[width=4.5cm]{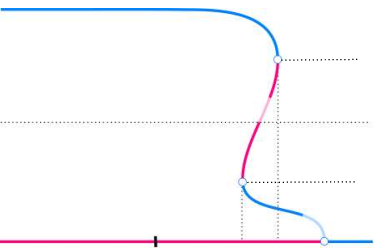}\vspace{0.5cm}
 \begin{picture}(0,0)
  \put(-358,-8){$0$}
  \put(-301,-8){$a_{H}$}
  \put(-405,48){\color{gray}with head}
  \put(-405,35){\color{gray}headless}
  \put(-220,-8){$0$}
  \put(-164,-8){$a_H$}
  \put(-80,-8){$0$}
  \put(-60,-8){$a^1_{sn}$}
  \put(-40,-8){$a^2_{sn}$}
  \put(-17,7){$a_H$}
  \put(-415,65){(a)}
  \put(-277,65){(b)}
  \put(-140,65){(c)}
  \put(-18,27){$x_1$}
  \put(-18,46){$-1$}
  \put(-18,70){$x_2$}
 \end{picture}
\end{center} 
 \caption{Representation of the width of limit cycles of system \eqref{systempal}-\eqref{vectorfield4z}, versus the parameter $a$, in the  supercritical case $m=-\sqrt{\varepsilon}$. The dotted line in all the panels corresponds with the width of the transitory canard, i.e. the limit between headless canard cycles and canard cycles with head. Moreover,  discoloured parts of the curves refers to the regions where the sign of the  functions $R_{3z}$ and $R_{4z}$ do not guarantee the stability of the limit cycle. In panels (a) and (b) we represent the cases where $k<1$ and $k=1$, respectively. In these cases, the limit cycle appearing after the supercritical bifurcation at $a_{H}$ exhibits a canard explosion. Panel (c) corresponds with the case $k>1$. Here, two saddle nodes of width $x_1<x_2$, take place at the values $a^1_{sn}<a^2_{sn}$ after the Hopf bifurcation at $a_{H}$.} \label{figteorneg}
\end{figure}

\begin{figure}[ht] 
\begin{center}
 \includegraphics[width=4.5cm]{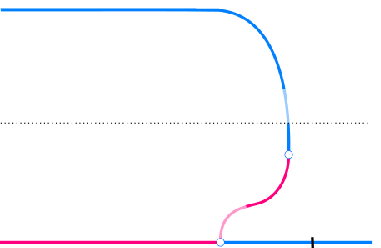}\quad \includegraphics[width=4.5cm]{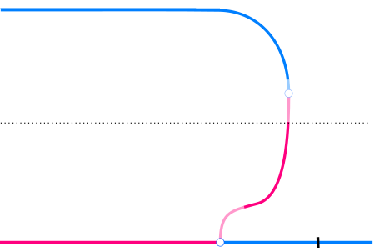}\quad \includegraphics[width=4.5cm]{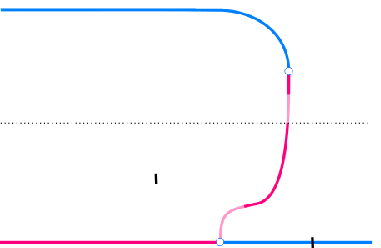}\vspace{0.5cm}
 \begin{picture}(0,0)
  \put(-338,-8){$a_{H}$}
  \put(-303,-8){$0$}
  \put(-405,48){\color{gray}with head}
  \put(-405,35){\color{gray}headless}
  \put(-200,-8){$a_H$}
  \put(-164,-8){$0$}
  \put(-27,-8){$0$}
  \put(-60,-8){$a_H$}
  \put(-415,65){(a)}
  \put(-277,65){(b)}
  \put(-139,65){(c)}
  \put(-20,48){-1}
 \end{picture}
\end{center} 
 \caption{Representation of the width of limit cycles of system \eqref{systempal}-\eqref{vectorfield4z}, versus the parameter $a$, in the  subcritical case $m=\sqrt{\varepsilon}$. The dotted line in all the panels corresponds with the transitory canard cycle with width $x=-1$. Moreover,  discolored parts of the curves refers to the regions where the sign of the functions $R_{3z}$ and $R_{4z}$ do not guarantee the stability of the limit cycle. In panels (a) and (b) we represent the cases where $k<1$ and $k=1$, respectively, and in panel (c) we represent the case $k>1$.} \label{figteorpos}
\end{figure}

In the next result we establish the stability of the canard limit cycles obtained in the previous theorem. To do this, we compute a piecewise smooth function 
\[
R(x)=\left\{
    \begin{array}{ll}
    R_{3z}(x) & x\in[-1,x_s),\\
    R_{4z}(x) & x\in (x_r,x_u),
    \end{array}
\right.
\]
see \eqref{cond1_3z} and \eqref{cond1_4z}, approximating the integral of the divergence along the limit cycle, $\Gamma_{x}$, and use the sign of this function to conclude the stability of $\Gamma_{x}$. Nevertheless, this approach does not produce accurate results if the canard limit cycle is close to the transitory canard, the one having width $x=-1$. 
The interval where the sign of $R(x)$  does not provide the stability of the canard limit cycles is given by $(x_u,-1)$, where 
\begin{equation}\label{def:xu}
x_u=-1+\lambda_{LL}^s(1+a).
\end{equation}
We organise the results into two theorems depending on whether the Hopf-like bifurcation is supercritical or subcritical. We illustrate the theorems with Figure \ref{figteorneg} and Figure \ref{figteorpos}.

\begin{theorem}\label{th:supercritico}
Set $\varepsilon>0$ small enough, $m=-\sqrt{\varepsilon}$, $x_0\in (x_r,x_u)\cup [-1,x_s)$ and $a=\hat{a}(k,\varepsilon,x_0;m)$. Let $\Gamma_{x_0}$ be the canard limit cycle of system \eqref{systempal}-\eqref{vectorfield4z} whose existence has been proved in Theorem \ref{th:exist_Gamma}. The following statements hold:
\begin{itemize}
\item [a)] For $k\leq 1$, the canard limit cycle $\Gamma_{x_0}$ is hyperbolic and stable. 
\item [b)] For $k>1,$ there exist exactly two values $x_1\in(-1,x_s)$ and $x_2\in(x_r,x_u)$ such that the canard limit cycle $\Gamma_{x_0}$ is hyperbolic and stable if $x_0\in(x_r,x_2)\cup(x_1,x_s),$ hyperbolic and unstable if $x_0\in(x_2,x_u)\cup(-1,x_1),$ and a saddle-node canard cycle if $x_0=x_1$ and $x_0=x_2.$
\end{itemize}
\end{theorem}

\begin{theorem}\label{th:subcritico}
Set $\varepsilon>0$ small enough, $m=\sqrt{\varepsilon}$, $x_0\in (x_r,x_u)\cup [-1,x_s)$ and $a=\hat{a}(k,\varepsilon,x_0;m)$. Let $\Gamma_{x_0}$ be the canard limit cycle of system \eqref{systempal}-\eqref{vectorfield4z} whose existence has been proved in Theorem \ref{th:exist_Gamma}. The following statements hold:
\begin{itemize}
\item [a)] For $k<1,$ there exists exactly one value $x_1\in(-1,x_s)$ such that $\Gamma_{x_0}$  is an hyperbolic limit cycle, if $x_0\in(x_r,x_u)\cup(-1,x_s)\setminus \{x_1\}$, and a saddle-node canard cycle, if $x_0=x_1$. Moreover, $\Gamma_{x_0}$ is stable if $x_0<x_1$ and unstable if $x_0>x_1.$
\item [b)] For $k=1,$ the canard limit cycle $\Gamma_{x_0}$ is hyperbolic and stable if $x_0\in(x_r,x_u)$ and hyperbolic and unstable if $x_0\in(-1,x_s).$
\item [c)] For $k>1,$ there exists exactly one value $x_2\in(x_r,x_u)$ such that $\Gamma_{x_0}$ is hyperbolic, if $x_0\in(x_r,x_u)\cup(-1,x_s)\setminus \{x_2\}$, and a saddle-node canard cycle, if $x_0=x_2$. Moreover, $\Gamma_{x_0}$ is stable if $x_0<x_2$ and unstable if $x_0>x_2.$
\end{itemize}
\end{theorem}

In the last main result we state that for every width between the smallest canard cycle and the relaxation oscillation cycle, that is for every $x_0\in(x_r,x_u)\cup [-1,x_s)$, there exist values of the parameters such that system \eqref{systempal}-\eqref{vectorfield4z} exhibits a saddle-node canard limit cycle $\Gamma_{x_0}$ of width $x_0$.

\begin{theorem}\label{th:supercritico_LP}
Consider system \eqref{systempal}-\eqref{vectorfield4z} with $m=-\sqrt{\varepsilon}$ or $m=\sqrt{\varepsilon}.$ For each $x_0\in(x_r,x_u)\cup(-1,x_s),$ there exists a value $\varepsilon_0$ and a function $k_{x_0}(\varepsilon)$ defined for $\varepsilon\in(0,\varepsilon_0),$ such that system \eqref{systempal}-\eqref{vectorfield4z} with parameters $k=k_{x_{0}}(\varepsilon)$ and $a=\hat{a}(k_{x_{0}}(\varepsilon),\varepsilon,x_0;m)$ exhibits the saddle-node canard $\Gamma_{x_{0}}$ whose existence has been stated in Theorem \ref{th:supercritico} for $m=-\sqrt{\varepsilon}$ and in Theorem \ref{th:subcritico} for $m=\sqrt{\varepsilon},$ respectively. 
\end{theorem}
%

\section{Proofs of the Main Results}\label{sec_proofs}
Let us begin by introducing some notation. For chosen parameters $\bm{\eta}=(a,k,m,\varepsilon)$, and a point $\mathbf{p} \in \mathbb{R}^2$, we denote by 
$$\varphi(t;\mathbf{p},\bm{\eta})=\left( x(t;\mathbf{p},\bm{\eta}), y(t;\mathbf{p},\bm{\eta})\right)$$ 
the solution of system \eqref{systempal} with initial condition $\varphi(0;\mathbf{p},\bm{\eta})=\mathbf{p}$. The coordinates of $\varphi(t;\mathbf{p},\bm{\eta})$ will be referred to as $x^i(t;\mathbf{p},\bm{\eta})$ and $y^i(t;\mathbf{p},\bm{\eta})$, with $i\in \{LL,L,C,R\}$, depending on the region where the solution belongs to, for that value of $t$.  


\subsection{Proof of Theorem \ref{th:connection}}\label{sec_theo1}

The existence of the maximal canard solution reduces to the existence of an orbit connecting points  $\mathbf{q}_1^R$ and $\mathbf{q}_0^L$, see Figure \ref{dib:geommetry}. The set of conditions characterizing this connection is given by the existence of ${\tau}_C>0$ and parameters $\bm{\eta}=(a,k,m,\varepsilon)$, with $0<\varepsilon\ll 1,$ $k>0$  and $|a|<\sqrt\varepsilon$, and such that the following conditions hold:
\begin{eqnarray}
 x^C(\tau_C; \mathbf{q}_1^R , \bm{\eta})&=&-\sqrt\varepsilon,\label{E1}\\
 y^C(\tau_C; \mathbf{q}_1^R, {\bm{\eta}})&=& -(m+\lambda_L^s)(\sqrt{\varepsilon}+a),\label{E2}\\
 x^C(s; \mathbf{q}_1^R, {\bm{\eta}})&\in&\left(-\sqrt\varepsilon,\sqrt\varepsilon\right)\; \mbox{ for all }\; s\in(0,{\tau}_C), \label{Ineq}
\end{eqnarray}
where $m=\pm \sqrt{\varepsilon}$. Next we proof the existence of a solution for system \eqref{E1}-\eqref{E2} and inequality \eqref{Ineq} in the supercritical case $m=-\sqrt{\varepsilon}$. The existence of a solution in the subcritical case, $m=\sqrt{\varepsilon}$, follows by completely analogous arguments. 

Let us define the following functions,
\begin{equation}\label{defFG}
\left\{\begin{array}{l}
F(\tau,a,k,\varepsilon)=x^C(\tau; \mathbf{q}_1^R , \bm{\eta})+\sqrt\varepsilon,\\ \noalign{\medskip}
G(\tau,a,k,\varepsilon)= y^C(\tau; \mathbf{q}_1^R, {\bm{\eta}})+(\lambda_L^s-\sqrt{\varepsilon})(\sqrt{\varepsilon}+a).
\end{array}\right.
\end{equation}
Thus, system (\ref{E1})-(\ref{E2}) is equivalent to system
\begin{equation}\label{sysFG}
\left\{\begin{array}{l}
F(\tau,a,k,\varepsilon)=0,\\
G(\tau,a,k,\varepsilon)=0,
\end{array}\right.
\end{equation}
and inequality (\ref{Ineq}) is equivalent to 
\begin{equation}
F(s,a,k,\varepsilon)\in\left(0,2\sqrt\varepsilon\right)\; \mbox{ for all }\; s\in(0,{\tau_C}). \label{IneqF}
\end{equation}

By integrating the linear system defined in the central band $C$ with initial condition ${\mathbf{q}}_1^R,$ we obtain the following explicit expression of the solution, 
\begin{equation}\label{solxy}
\left\{\begin{array}{ll}
x^C(\tau;\mathbf{q}_1^R,\bm{\eta})&=
 \displaystyle \frac {e^{\frac{\sqrt{\varepsilon}}2\tau}(\sqrt{\varepsilon}-a)}{\sqrt{3}}
  \left(
    \left(
      \frac{2\lambda_R^s}{\sqrt{\varepsilon}}-1
    \right) \sin\left(\frac{\sqrt{3\varepsilon}}{2} \tau\right)
    +\sqrt{3}\cos\left(\frac{\sqrt{3\varepsilon}}{2} \tau\right)
  \right)+a \\  \noalign{\medskip} 
y^C(\tau;\mathbf{q}_1^R,\bm{\eta})&=
  \displaystyle \frac {e^{\frac{\sqrt{\varepsilon}}2\tau}(\sqrt{\varepsilon}-a)}{\sqrt{3}}
  \left((\lambda_R^s-\sqrt{\varepsilon})\sqrt{3} \cos\left(\frac{\sqrt{3\varepsilon}}{2} \tau\right)
 -(\lambda_R^s+\sqrt{\varepsilon}) \sin\left(\frac{\sqrt{3\varepsilon}}{2} \tau\right) \right).
\end{array}\right.
\end{equation}

Consider the rescaling of parameters $\tau$ and $a,$ as follows :
\begin{equation}\label{rescaling}
\quad \tau=\bar\tau/\sqrt\varepsilon,\quad a=\bar a\sqrt{\varepsilon}.
\end{equation}
It follows from the definition of $F$ (see (\ref{defFG})) and formula for $x^C$ (see (\ref{solxy})) that $F=\sqrt\varepsilon\bar F$, where $\bar F$ is given by  
\begin{equation}\label{eq-hatf}
\bar F(\bar{\tau},\bar a,\varepsilon)=\sqrt{3} (\bar{a}-1) e^{\bar{\tau}/2} \left(\left(2
   \sqrt{\varepsilon}+1\right) \sin \left(\frac{\sqrt{3}
   \bar{\tau}}{2}\right)-\sqrt{3} \cos \left(\frac{\sqrt{3}
   \bar{\tau}}{2}\right)\right)+3 \bar{a}+3+O(\varepsilon).
\end{equation}
In the same way, from the definition of $G$ (see (\ref{defFG})) and expressions for $x^C$ and $y^C$ (see (\ref{solxy})), it holds that $G=\varepsilon\bar G$, where $\bar G$ is given by  
\begin{equation*}
\bar G(\bar\tau,\bar a,k,\varepsilon)=\frac{(\bar{a}+1)
   \left(\sqrt{\varepsilon}-k\right)}{k}+\frac{(\bar{a}-1)
   e^{\bar{\tau}/2} \left(\sqrt{3}
   \left(\sqrt{\varepsilon}+1\right) \cos \left(\frac{\sqrt{3}
   \bar{\tau}}{2}\right)-\left(\sqrt{\varepsilon}-1\right) \sin
   \left(\frac{\sqrt{3} \bar{\tau}}{2}\right)\right)}{\sqrt{3}}+O(\varepsilon).
\end{equation*}
Hence, solving $F=G=0$ for $\varepsilon> 0$ is equivalent to solving $\bar F=\bar G=0$.
For $\varepsilon=0,$ we have
\begin{equation}\begin{array}{ll}\label{computationFG}
\bar F&=\sqrt{3} (\bar{a}-1) e^{\bar{\tau}/2} \left(\sin \left(\frac{\sqrt{3}
   \bar{\tau}}{2}\right)-\sqrt{3} \cos \left(\frac{\sqrt{3}
   \bar{\tau}}{2}\right)\right)+3 \bar{a}+3,\\ \noalign{\smallskip}
\bar G&=-(\bar{a}+1)+\displaystyle\frac{(\bar{a}-1)
   e^{\bar{\tau}/2} \left(\sqrt{3}
  \cos \left(\frac{\sqrt{3}
   \bar{\tau}}{2}\right)+\sin
   \left(\frac{\sqrt{3} \bar{\tau}}{2}\right)\right)}{\sqrt{3}}.
\end{array}\end{equation}
For $\bar a=(e^{\frac{\pi}{\sqrt{3}}}-1)/(e^{\frac{\pi}{\sqrt{3}}}+1), \bar\tau=2\pi/\sqrt{3}$ equations $\bar F=0,\bar G=0$ hold. To apply the implicit function theorem, 
it is necessary to prove that $\det(J(2\pi/\sqrt{3},(e^{\frac{\pi}{\sqrt{3}}}-1)/(e^{\frac{\pi}{\sqrt{3}}}+1),k,0))\neq 0,$ where
\begin{equation}J(\bar \tau,\bar a,k,\varepsilon)=
\left(
\begin{array}{cc}
  \displaystyle\frac{\partial \bar F}{\partial \bar \tau}(\bar \tau,\bar a,k,\varepsilon)&  \displaystyle\frac{\partial \bar F}{\partial \bar  a}(\bar \tau,\bar a,k,\varepsilon)\medskip\\
    \displaystyle\frac{\partial \bar G}{\partial \bar \tau}(\bar \tau,\bar a,k,\varepsilon)&  \displaystyle\frac{\partial \bar G}{\partial \bar a}(\bar \tau,\bar a,k,\varepsilon)\\
  \end{array}
\right).
\end{equation}
It is easy to see that the partial derivatives take the following values at the point  $(2\pi/\sqrt{3},(e^{\frac{\pi}{\sqrt{3}}}-1)/(e^{\frac{\pi}{\sqrt{3}}}+1),k,0),$
\[
\frac{\partial \bar F}{\partial \bar  \tau}=0,\quad \frac{\partial \bar  F}{\partial \bar  a}=3(1+e^{\frac{\pi}{\sqrt{3}}}),
\]
\[
\frac{\partial \bar G}{\partial \bar \tau}=\frac{2e^{\frac{\pi}{\sqrt{3}}}}{1+e^{\frac{\pi}{\sqrt{3}}}},\quad \frac{\partial \bar G}{\partial \bar a}=-1-e^{\frac{\pi}{\sqrt{3}}},
\]
and then, $\det(J(2\pi/\sqrt{3},(e^{\frac{\pi}{\sqrt{3}}}-1)/(e^{\frac{\pi}{\sqrt{3}}}+1),k,0))=-6e^{2\frac{\pi}{\sqrt{3}}}\neq 0.$
Thus, from the implicit function theorem we conclude that there exist a value $\varepsilon_0>0$ and functions $\tilde{a}(k,\varepsilon;m)$ and $\tau_C(k,\varepsilon;m)$, analytic and smooth as a function of $(k,\sqrt{\varepsilon})$, respectively, defined in the open set $U=(0,+\infty)\times (0,\varepsilon_0)$ and such that equations (\ref{E1})-(\ref{E2}) has the solution $a=\tilde{a}(k,\varepsilon;m)$ and $\tau=\tau_C(k,\varepsilon;m)$. Moreover, the lower order terms in $\sqrt{\varepsilon}$ of such a solution coincides with those in the expression (\ref{aepmenos})-(\ref{tauepmenos}). The remainder terms in this approximated solution can be obtained by the method of the undetermined coeficients. 

Finally, since the angle travelled by the solution from $\mathbf{q}_1^R$ to $\mathbf{q}_0^L$ satisfies that  $\beta \tau_C(k,\varepsilon;m)<\pi$ and $|\tilde{a}(k,\varepsilon;m)|<\sqrt{\varepsilon}$, it follows that inequality (\ref{Ineq}) is fulfilled. Therefore, we conclude that for $a=\tilde{a}(k,\varepsilon;m)$ system \eqref{systempal}-\eqref{vectorfield4z} exhibits an orbit connecting the slow manifolds with time of flight equal to $\tau_C(k,\varepsilon;m)>0$. 

\subsection{Proof of Theorem \ref{th:exist_Gamma}}
We perform the proof for the headless canard limit cycles. The proof for cycles with head follows in a similar way. Consider a point $(x_0,f(x_0))$ with  $x_0\in(-1,{x}_s)$, and the orbit $\Gamma_{x_0}$ through this point, see Figure \ref{dib:geommetry}. Since the point is contained in the positive invariant region $\mathcal{R}$, see Lemma \ref{lem:inv_reg}, then by integrating both in forward and backward time, the orbit targets the switching line $\{x=-\sqrt{\varepsilon}\}$ at two points respectively
\[
\begin{pmatrix}
-\sqrt{\varepsilon}\\h
\end{pmatrix},\quad\quad
\mathbf{p}_0 = \mathbf{q}_0^L+
\begin{pmatrix}
0 \\ h e^{-\frac {k h}{\varepsilon(\sqrt{\varepsilon}-a)}}
\end{pmatrix}.
\]
The expression of $\mathbf{p}_0$ follows from Lemma \ref{lem:Pmaps}, since $h=\Phi_{3z}(x_0) > (-m+\varepsilon^{-\nu}\lambda_L^s)(\sqrt{\varepsilon}+a)$ provided that $x_0\in (-1,x_s)$, see Lemma \ref{lem:dom}.
Moreover, $\Gamma_{x_0}$ also intersects in forward time the switching line $\{x=\sqrt{\varepsilon}\}$ at the two points 
\[
\begin{pmatrix}
\sqrt{\varepsilon}\\h
\end{pmatrix},\quad\quad
\mathbf{p}_1=\mathbf{q}_1^R+
\begin{pmatrix}
0 \\ h e^{-\frac {h}{\varepsilon(\sqrt{\varepsilon}-a)}}
\end{pmatrix}.
\]
The conditions on the points $\mathbf{p}_0$ and $\mathbf{p}_1$ in order to be $\Gamma_{x_0}$ a limit cycle is the existence of values for the parameters $\bm{\eta}=(a,k,m,\varepsilon)$ such that the solution of the linear system in the central band $C$ with initial condition at $\mathbf{p}_1$ targets first time the switching line $\{x=-\sqrt{\varepsilon}\}$ at the point $\mathbf{p}_0$, i.e., $e^{\tau A_c}(\mathbf{p}_1-\mathbf{e}_c) + \mathbf{e}_c-\mathbf{p}_0 = \mathbf{0}$ with $\beta \tau<\pi$. By substituting the values of $\mathbf{p}_0$ and $\mathbf{p}_1$, previous equation writes as
\begin{equation}\label{aux0:lem_43}
e^{\tau A_c}(\mathbf{q}_1^R -\mathbf{e}_c) + \mathbf{e}_c-\mathbf{q}_0 + 
e^{\tau A_c} \begin{pmatrix} 0 \\ h e^{-\frac {h}{\varepsilon(\sqrt{\varepsilon}-a)}} \end{pmatrix}
-\begin{pmatrix}
0 \\ h e^{-\frac {k h}{\varepsilon(\sqrt{\varepsilon}-a)}}
\end{pmatrix}
= \mathbf{0}.
\end{equation}
Considering the change of variables $\tau=\frac{\tilde{\tau}}{\sqrt{\varepsilon}}$ and $a=\hat{a}\sqrt{\varepsilon}$ and multiplying previous equations by $\frac {1}{\sqrt{\varepsilon}}$ and $\frac 1{\varepsilon}$, respectively, we obtain that expression \eqref{aux0:lem_43}  writes as $E_p(\tilde{\tau},\hat{a},k,m,\varepsilon)=\mathbf{0}$, where
\begin{equation}\label{aux:lem_43}
E_p(\tilde{\tau},\hat{a},k,m,\varepsilon,h)= E_q(\tilde{\tau},\hat{a},k,m,\varepsilon) +
\begin{pmatrix}
\frac 1{\sqrt{\varepsilon}} & 0 \\
0 & \frac 1{\varepsilon} 
\end{pmatrix}
\left(
e^{\tilde{\tau} A_c} \begin{pmatrix} 0 \\ h e^{-\frac {h}{\varepsilon(\sqrt{\varepsilon}-\hat{a})}} \end{pmatrix}
-\begin{pmatrix}
0 \\ h e^{-\frac {k h}{\varepsilon(\sqrt{\varepsilon}-\hat{a})}}
\end{pmatrix}
\right),
\end{equation}
and $E_q(\tilde{\tau},\hat{a},k,m,\varepsilon)=e^{\tilde{\tau} A_c}(\mathbf{q}_1^R -\mathbf{e}_c) + \mathbf{e}_c-\mathbf{q}_0$ are the equations given in (\ref{computationFG}) and establishing the connection between $\mathbf{q}_1^R$ and  $\mathbf{q}_0^L$. Then, the Jacobian matrix respect to the variables $\tau$ and $a$ satisfies
\begin{align*}
\left. D_{\tilde{\tau},\hat{a}}E_p\right|_{(\tilde{\tau},\hat{a},k,m,\varepsilon,h)}&= \left.D_{\tilde{\tau},\hat{a}}E_q\right|_{(\tilde{\tau},\hat{a},k,m,\varepsilon)}\\ 
&+ \begin{pmatrix}
\frac 1{\sqrt{\varepsilon}} & 0 \\
0 & \frac 1{\varepsilon} 
\end{pmatrix}
\left(
e^{\tilde{\tau} A_c} 
\begin{pmatrix} 
   h e^{-\frac {h}{\varepsilon(\sqrt{\varepsilon}-\hat{a})}} & 0 \\
   0 & -\frac{h^2 e^{-\frac {h}{\varepsilon(\sqrt{\varepsilon}-\hat{a})}}}{\varepsilon(\sqrt{\varepsilon}-\hat{a})^2}
\end{pmatrix} -
\begin{pmatrix} 
   0 \\ -\frac{kh^2 e^{-\frac {k h}{\varepsilon(\sqrt{\varepsilon}-\hat{a})}}}
   {\varepsilon(\sqrt{\varepsilon}-\hat{a})^2}    
\end{pmatrix}
\right).
\end{align*}

The trace $t_c=-\sqrt{\varepsilon}$ and the determinant $d_c=\varepsilon$ of the matrix $A_C$ tend  to zero when $\varepsilon$ tends to zero. Hence, for every $h$ with 
\begin{equation*}\label{aux2:lem43}
    h>O(\varepsilon^{3/2}),
\end{equation*}
we conclude that  $E_p(\tilde{\tau},\hat{a},k,0,0,h)= E_q(\tilde{\tau},\hat{a},k,0,0)$ and 
$\left. D_{\tilde{\tau},\hat{a}}E_p\right|_{(\tilde{\tau},\hat{a},k,0,0,h)}= \left.D_{\tilde{\tau},\hat{a}}E_q\right|_{(\tilde{\tau},\hat{a},k,0,0)}$. 

From the proof of Theorem \ref{th:connection}, it follows that
\begin{align*}
E_p(2\pi/\sqrt{3},(e^{\frac{\pi}{\sqrt{3}}}-1)/(e^{\frac{\pi}{\sqrt{3}}}+1),k,0,0,h)&=\mathbf{0},\\
\det\left(\left.D_{\tilde{\tau},\hat{a}}E_p\right|_{(2\pi/\sqrt{3},(e^{\frac{\pi}{\sqrt{3}}}-1)/(e^{\frac{\pi}{\sqrt{3}}}+1),k,0,0,h)}\right)&= -6e^{\frac{2\pi}{\sqrt{3}}}.
\end{align*}
Thus, we can apply the Implicit Function Theorem to the set of equations $E_p(\tilde{\tau},\hat{a},k,m,\varepsilon,h)=\mathbf{0}$, which establishes the connection between $\mathbf{p}_{1}$  and $\mathbf{p}_{0}$, and conclude the existence of the functions $\tilde{\tau}=\tilde{\tau}_C(k,\varepsilon,h;m)$ and $\hat{a}=\hat{a}(k,\varepsilon,h;m)$ satisfying 
\[
E_p(\tilde{\tau}_C(k,\varepsilon,h;m),\hat{a}(k,\varepsilon,h;m),k,m,\varepsilon,h)=\mathbf{0}.
\]

Furthermore, it follows that
\begin{align*}
E_p(\tilde{\tau}_C(k,\varepsilon,h;m),\hat{a}(k,\varepsilon,h;m),k,m,\varepsilon,h) & -
E_p(\tau_C(k,\varepsilon;m),\tilde{a}(k,\varepsilon;m),k,m,\varepsilon,h) = \\
\tilde{E}_p(\tilde{\tau}_C,\hat{a},\tau_C,\tilde{a},k,m,\varepsilon,h)
&\begin{pmatrix}
\tilde{\tau}_C(k,\varepsilon,h;m) - \tau_C(k,\varepsilon;m)\\
\hat{a}(k,\varepsilon,h;m)-\tilde{a}(k,\varepsilon;m)
\end{pmatrix},
\end{align*}
where $\tilde{E}_p(\tilde{\tau}_C,\hat{a},\tau_C,\tilde{a},k,m,\varepsilon,h)=\int_0^1
D_{\tilde{\tau},\hat{a}}E_p(s\tilde{\tau}_C+(1-s)\tau_C,s\hat{a}+(1-s)\tilde{a},k,m,\varepsilon,h)\, ds$.

Applying now equation \eqref{aux:lem_43}, we conclude that for $\varepsilon$ small enough 
\[
\begin{pmatrix}
\tilde{\tau}_C(k,\varepsilon,h;m) - \tau_C(k,\varepsilon;m)\\
\hat{a}(k,\varepsilon,h;m)-\tilde{a}(k,\varepsilon;m)
\end{pmatrix}
\approx
\frac {6e^{\frac{2\pi}{\sqrt{3}}}}{ \varepsilon} 
\left(
e^{\bar{\tau} Ac}
\begin{pmatrix}
0 \\ he^{-\frac{h}{\varepsilon^{3/2}}}
\end{pmatrix}-
\begin{pmatrix}
0 \\ he^{-\frac{k h}{\varepsilon^{3/2}}}
\end{pmatrix}
\right),
\]
where $\bar{\tau}=(e^{\frac{\pi}{\sqrt{3}}}-1)/(e^{\frac{\pi}{\sqrt{3}}}+1)$. From here we conclude the proof.


\subsection{Proof of theorems \ref{th:supercritico} and \ref{th:subcritico}}\label{sec:3z}
This subsection is devoted to the proof of theorems \ref{th:supercritico} and \ref{th:subcritico}. The outline of the proof is: First, we study the hyperbolicity/non-hyperbolicity of headless canard cycles. Second, we analyze the  hyperbolicity/non-hyperbolicity of canard cycles with head. Finally, we prove the correspondence of non-hyperbolic canard cycles to saddle-node bifurcations.

\subsubsection{Hyperbolicity/non-hyperbolicity of headless canard cycles}

Let us begin by defining the Poincaré map in neighborhood of periodic orbits visiting zones $L$, $C$ and $R$.
Consider a point $\mathbf{p}_0=(\sqrt{\varepsilon},y_0)$  in the switching line $\{x=\sqrt{\varepsilon}\}$ and located between $\mathbf{p}_R$ and $\mathbf{q}_{R}^1$, see Figure \ref{dib:geommetry}. From expressions \eqref{q1R} and \eqref{pLpLL} it follows that $(m+\lambda_R^s)(\sqrt{\varepsilon}-a)<y_0<m(\sqrt{\varepsilon}-a)$. Assume now that there exists a time of flight $\tau_{Cd}>0$ such that 
$x^C(\tau_{Cd};\mathbf{p}_0,\bm{\eta})=-\sqrt{\varepsilon}$ and  $x^C(s;\mathbf{p}_0,\bm{\eta})\in(-\sqrt{\varepsilon},\sqrt{\varepsilon})$ for all $s\in(0,\tau_{Cd})$, where $x^C(s;\mathbf{p}_0,\bm{\eta})$ is the first coordinate of the solution through $\mathbf{p}_0$ reduced to the central region $C$, see \eqref{E1}.
In such a case, we can define the Poincar\'e half-map between the switching lines $\{x=\sqrt{\varepsilon}\}$ and $\{x=-\sqrt{\varepsilon}\}$ at the point $y_0$ as $\Pi_{C_d}(y_{0},{\bm{\eta}})= {y^C(\tau_{Cd};\mathbf{p}_0,{\bm{\eta}})}.$
Similarly, we can define the Poincar\'e half-map between the switching lines $\{x=-\sqrt{\varepsilon}\}$ and $\{x=\sqrt\varepsilon\}$ at a point $y_2>-m(\sqrt{\varepsilon}+a)$ as
$\Pi_{C_u}(y_2,{\bm{\eta}})=y^C\left(\tau_{C_u};\mathbf{p}_2,{\bm{\eta}}\right)$, where $\tau_{C_u}>0$ is the time of flight and $\mathbf{p}_2=(-\sqrt{\varepsilon},y_2).$ 

Consider now a point $\mathbf{p}_1=(-\sqrt{\varepsilon},y_1)$ in the switching line $\{x=-\sqrt{\varepsilon}\}$ and located between $\mathbf{q}_L^0$ and $\mathbf{p}_L$, that is $y_1\in\left(-(m+\lambda_L^s)(\sqrt{\varepsilon}+a),-m(\sqrt{\varepsilon}+a)\right)$.  Assume that there exists a time of flight $\tau_{L}>0$ such that $x^L(\tau_{L};\mathbf{p}_1,\bm{\eta})=-\sqrt{\varepsilon}$ and  $x^L(s;\mathbf{p}_1,\bm{\eta})\in(-1,-\sqrt{\varepsilon})$ for all $s\in(0,\tau_L).$ Here $x^L(\tau_{L};\mathbf{p}_1,\bm{\eta})$ is the first coordinate of the solution through $\mathbf{p}_1$ and reduced to the region $L$. In such a case, we can define the Poincaré half-map between the switching line $\{x=-\sqrt{\varepsilon}\}$ and itself at the point $y_1$ as $\Pi_L(y_1,{\bm{\eta}})=y^L\left(\tau_{L};\mathbf{p}_1,{\bm{\eta}}\right)$.
Similarly, we can define the Poincar\'e half-map between the switching line $\{x=\sqrt{\varepsilon}\}$ and itself at the point $y_3>m(\sqrt{\varepsilon}-a)$ as
$\Pi_R(y_3,{\bm{\eta}})=y^R\left(\tau_{R};\mathbf{p}_3,{\bm{\eta}}\right)$, where $\tau_{R}>0$ is the time of flight and $\mathbf{p}_3=(\sqrt{\varepsilon},y_3).$
Expressions for $\Pi_L^{-1}$ and $\Pi_R$ can be found in Lemma \ref{lem:Pmaps}.

At this point, the Poincar\'{e} map for an orbit of system \eqref{systempal} 
visiting zones $L$, $C$ and $R$ can be defined.
\begin{definition}\label{defipoincareydesplazamiento}
The Poincar\'e map $\Pi$ in the neighborhood of an orbit $\Gamma_{x_0}$ of system \eqref{systempal} 
visiting zones $L$, $C$ and $R$  is defined as
\begin{equation*}
\Pi(y_0,\bm{\eta})=\Pi_R(\Pi_{C_u}(\Pi_L(\Pi_{C_d}(y_0,\bm{\eta}),\bm{\eta}),\bm{\eta}),\bm{\eta}),
\end{equation*}
provided the composition of Poincar\'e half-maps is possible, where $y_0=\Pi_{Cd}^{-1}(\Pi_{L}^{-1}(\Phi_{3z}(x_0),{\bm \eta}),{\bm \eta})$.
\end{definition}

For $\varepsilon$ fixed and small enough, suppose the existence of a headless canard limit cycle $\Gamma_{x_0}$, see Figure \ref{dib:geommetry}, obtained under the parameter relation $a=\hat{a}(k,\varepsilon,x_0; m)$ given in Theorem \ref{th:exist_Gamma}. 
The cycle  $\Gamma_{x_0}$ corresponds to the fixed point of the Poincar\'e map $\Pi(y_0,{\bm  \eta}),$ where $y_0=\Pi_{Cd}^{-1}(\Pi_{L}^{-1}(\Phi_{3z}(x_0),{\bm \eta}),{\bm \eta})$.

To take into account the non-hyperbolicity of the canard cycle $\Gamma_{x_0}$, we consider the derivative of the Poincar\'e map, which corresponds to the exponential of the integral of the divergence of the system along $\Gamma_{x_0}$, see  \cite{chicone}. In the particular case of PWL systems, the integral of the divergence can be explicitly computed as the sum of the products of the traces and the time of flight of $\Gamma_{x_0}$ in each region of linearity, see \cite{freire1}. 

\begin{figure}[h]
    \centering
    \includegraphics{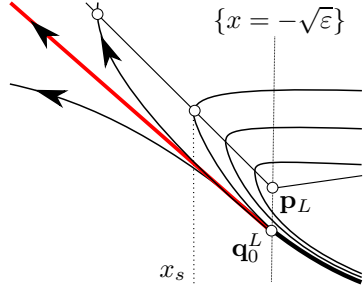}
    \begin{picture}(0,0)
    \put(-80,4){$x_s$}
    \put(-53,12){$\mathbf{q}_0^L$}
    \put(-35,30){$\mathbf{p}_L$}
    \put(-60,98){$\{x=-\sqrt{\varepsilon}\}$}
    \end{picture}
    \caption{Zoom of the flow in a neighbourhood of the contact point $\mathbf{p}_L$. Orbits having width in $x_s<x<-\sqrt{\varepsilon}$ do not pass exponentially close of $\mathbf{q}_0^L$. Therefore, the time of flight $\tau_L$ can not be computed as it is done in Lemma \ref{lem:times}.}
    \label{fig:Zoom1}
\end{figure}

Let $\tau_L$ and $\tau_R$ be the time of flight of $\Gamma_{x_0}$ along the regions $L$ and $R$, respectively, and let $\tau_C=\tau_C(k,\varepsilon;m)$ be the time of flight from $\mathbf{q}_1^R$ to $\mathbf{q}_0^L$ obtained in Theorem \ref{th:connection}. From Lemma \ref{lem:dom}, when $x_0\in[-1,x_s)$ it follows that $\Gamma_{x_0}$ intersect the switching lines $\{x=-\sqrt{\varepsilon}\}$ and $\{x=\sqrt{\varepsilon}\}$ exponentially close to $\mathbf{q}_0^L$ and $\mathbf{q}_1^R$, respectively. Therefore, the values of $\tau_L$ and $\tau_R$ can be approximated by the time of flight of the orbit from $\mathbf{q}_0^L$ to $(-\sqrt{\varepsilon},h)$ and from $(\sqrt{\varepsilon},h)$ to $\mathbf{q}_1^R$, respectively. Hence, $\tau_L=\tau_L(h)$ and $\tau_R=\tau_R(h)$ are the ones computed in Lemma \ref{lem:times}. Notice that when $x_s<x_0<-\sqrt{\varepsilon}$, we can not assure that $\Gamma_{x_0}$ intersect neither $\{x=-\sqrt{\varepsilon}\}$ nor $\{x=\sqrt{\varepsilon}\}$ exponentially close to $\mathbf{q}_0^L$ and $\mathbf{q}_1^R$, respectively, see Figure \ref{fig:Zoom1}.
In such case, expressions in Lemma \ref{lem:times} are not good approximations to $\tau_L$ and $\tau_R$, respectively. Therefore we have eliminated the interval $(x_s,-\sqrt{\varepsilon})$ from the scope of the Theorem \ref{th:supercritico} and of the Theorem \ref{th:subcritico}. 

Finally, the time of flight from $(-\sqrt{\varepsilon},h)$ to $(\sqrt{\varepsilon},h)$ (dashed arc on $\Gamma_{x_0}$ in Figure \ref{dib:geommetry}) is order one in $\varepsilon$, so we can consider it as zero.

Then, the derivative of the Poincaré map $\Pi(y_0,{\bm \eta})$, where $y_0=\Pi_{Cd}^{-1}(\Pi_{L}^{-1}(h,{\bm \eta}),{\bm \eta})$ and $h=\Phi_{3z}(x_0)$, can be approximated by 
\begin{equation}\label{derprim}
\frac{\partial \Pi}{\partial y}(y_0,\bar{\bm{\eta}})\approx e^{t_L \tau_L-m\tau_C+t_R\tau_R}. 
\end{equation}

A necessary condition on the limit cycle $\Gamma_{x_0}$ to be non-hyperbolic is that 
$\frac{\partial \Pi}{\partial y}(y_0,\bar{\bm{\eta}})=1$ with $y_0=\Pi_{Cd}^{-1}(\Pi_{L}^{-1}(\Phi_{3z}(x_0),{\bm \eta}),{\bm \eta})$. In terms of the right side of \eqref{derprim} this condition writes as $R_{3z}(h,k,\varepsilon;m)=e^{t_L \tau_L+t_R\tau_R}-e^{m\tau_C}=0.$ 
By using the expression of $\tau_L$, $\tau_R$ and $\tau_C$ previously computed it follows that 
\begin{equation}\label{cond1_3z}
  R_{3z}(h,k,\varepsilon;m)=
  \left(1 + \frac {h+(m+\lambda_L^s)(\sqrt{\varepsilon}+a)}{(\lambda_L^q-\lambda_L^s)(\sqrt{\varepsilon}+a)}  \right)^{\frac k{\lambda_L^s}} 
  \left(1 +  \frac {(m+\lambda_R^s)(\sqrt{\varepsilon}-a)-h}{(\lambda_R^q-\lambda_R^s)(\sqrt{\varepsilon}-a)} \right)^{\frac 1{\lambda_R^s}}-
  e^{m\tau_C}.
\end{equation}
Notice that the function $R_{3z}$ writes in terms of the height $h$ of the cycle $\Gamma_{x_0}$, with $h=\Phi_{3z}(x_0)$. In the next result we compute the stability of a canard cycle through the sign of the function $R_{3z}$.


\begin{proposition}\label{prop:R3z}
For $\varepsilon$ fixed and small enough, there exists $0<\delta\ll 1$ such that for $x_0\in[-1,x_s)$ and $h=\Phi_{3z}(x_0)$: 
\begin{itemize}
 \item [a)] if $R_{3z}(h,k,\sqrt{\varepsilon}; m)<-\delta$, then $\Gamma_{x_0}$ is a hyperbolic headless stable canard cycle;
 \item [b)] if $h$ is a simple root of $R_{3z}(h,k,\sqrt{\varepsilon}; m)$, then in a neighborhood of $\Gamma_{x_0}$ there is a nonhyperbolic headless canard cycle;
 \item [c)] if $R_{3z}(h,k,\sqrt{\varepsilon}; m)>\delta$, then $\Gamma_{x_0}$ is a hyperbolic headless unstable canard cycle.
\end{itemize}
\end{proposition}
\begin{proof}
The proposition is a straight forward consequence of the definition of $R_{3z}(h,k,\varepsilon;m)$ and the equation \eqref{derprim}.
\end{proof}

By fixing parameters $k$ and $\varepsilon$ we next describe the qualitative behavior of $R_{3z}(h,k,\varepsilon;m)$ as a function of $h$, both in the supercritical case,  $m=-\sqrt{\varepsilon}$, and in the subcritical case, $m=\sqrt{\varepsilon}$. Even when the  domain of definition of $R_{3z}$, as a function of $h$, is greater, we only consider the reduction of $R_{3z}$ to the interval $(h_s,h_M]$ where  $h_s=\Phi_{3z}(x_s)$ and $h_M=\Phi_{3z}(-1)$, see Lemma \ref{lem:dom}. We pay special attention to the existence of simple zeros of $R_{3z}$.

\begin{figure}[ht]
\begin{center}
\includegraphics{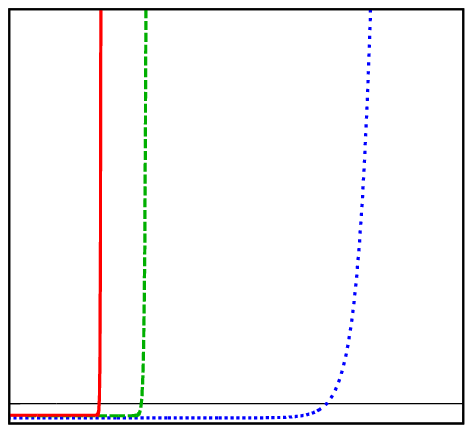}\quad\quad\includegraphics{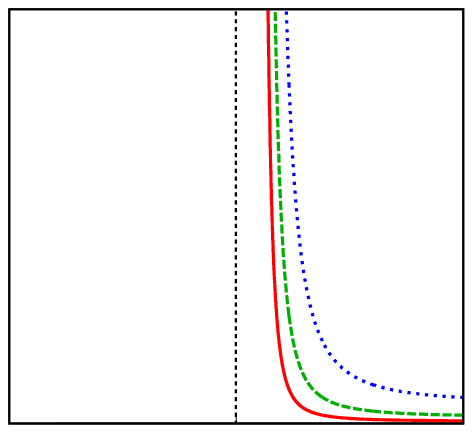}\\
\includegraphics{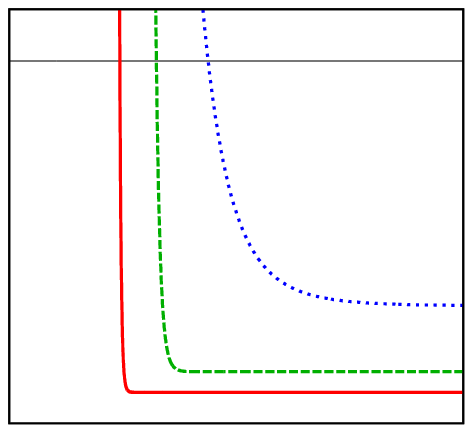}\quad\quad\includegraphics{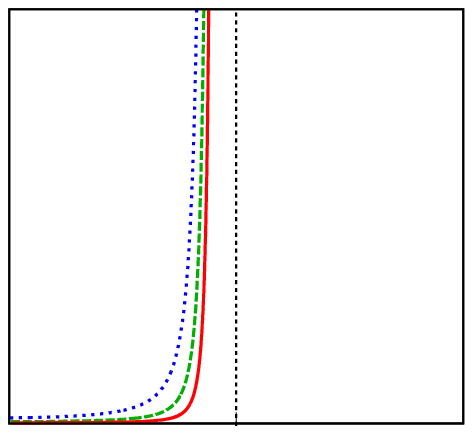}
\begin{picture}(0,0)
\put(-360,275){(a)}
\put(-236,254){\color{white}\circle*{15}}
\put(-280,250){{\small $R_{3z}(h,k,\sqrt{\varepsilon};m)$}}
\put(-220,147){$h$}
\put(-352,165){$0$}
\put(-248,165){\circle*{4}}
\put(-302,165){\circle*{4}}
\put(-314,165){\circle*{4}}
\put(-255,170){$h^*$}
\put(-170,275){(b)} 
\put(-86,147){$1$}
\put(-30,147){$k$}
\put(-165,245){$h^*$}
\put(-60,147){$2.5$}
\put(-55,160){\line(0,1){25}}
\put(-55,162){\circle*{4}}
\put(-55,166){\circle*{4}}
\put(-55,178){\circle*{4}}
\put(-65,204){\color{white}\circle*{15}}
\put(-69,200){{\small $h^*(k,\sqrt{\varepsilon};m)$}}
\put(-360,130){(c)}
\put(-352,115){$0$}
\put(-278,85){\color{white}\circle*{15}}
\put(-290,82){{\small $R_{3z}(h,k,\sqrt{\varepsilon};m)$}}
\put(-220,4){$h$}
\put(-170,130){(d)}
\put(-86,4){$1$}
\put(-30,4){$k$}
\put(-165,100){$h^*$}
\put(-104,64){\color{white}\circle*{15}}
\put(-145,61){{\small $h^*(k,\sqrt{\varepsilon};m)$}}
\end{picture}
\end{center}
\caption{Graphical representation of the functions $R_{3z}(h,k,\sqrt{\varepsilon};m)$ and $h^*(k,{\varepsilon};m)$ analysed in Proposition \ref{th:nhcond3z} for different values of $\varepsilon$. In particular, blue/dotted curves correspond with $\varepsilon=0.05$, green/dashed curves with  $\varepsilon=0.01$ and red/solid curves with $\varepsilon=0.005$. First column contains the graphs of the function $R_{3z}(h,k,\varepsilon;m)$ as a function of $h$ for the previous values of $\varepsilon$: panel (a) when $k=2.5$ and $m=-\sqrt{\varepsilon}$ and panel (c) when $k=0.75$ and $m=\sqrt{\varepsilon}$. Second column contains the graphs of the function $h^*(k,\varepsilon;m)$ as a function of $k$ for the previous values of $\varepsilon$: panel (b) when $m=-\sqrt{\varepsilon}$ and panel (c) when $m=\sqrt{\varepsilon}$.}\label{dib:FunsSyh}
\end{figure}

\begin{proposition}\label{th:nhcond3z}
Fixed $\varepsilon$ small enough, we consider the  function $R_{3z}(h,k,\varepsilon;m)$ defined in 
\eqref{cond1_3z}.
\begin{itemize}
  \item [a)] Under the supercritical condition $m=-\sqrt{\varepsilon}$ we obtain that:
  \begin{itemize} 
    \item [a-1)] for $k\leq 1$, then $R_{3z}(h,k,\varepsilon;m)<0$ when $h\in(h_s,h_M]$, and 
    \item [a-2)] for $k>1$, the function $R_{3z}(h,k,\varepsilon;m)$ behaves as it is represented in Figure \ref{dib:FunsSyh}(a). More specifically, 
      \begin{itemize}
	\item [a-2-1)] $\lim_{h\searrow h_s}R_{3z}(h,k,\varepsilon;m)<0$ and $ R_{3z}(h_M,k,\sqrt{\varepsilon};m)>0$, 
	\item [a-2-2)] let $h^*\in(h_s,h_M]$ be a zero of $R_{3z}(h,k,\varepsilon;m)$, then $ \left. \frac {\partial R_{3z}}{\partial h} \right|_{(h^*,k,\sqrt{\varepsilon};m)} >0$, and
	\item [a-2-3)] denoting by $h^*(k,\varepsilon;m)$ the unique positive zero of $R_{3z}(h,k,\varepsilon;m)$ in $(h_s,h_M]$, then
	\[
	  h^*(k,\varepsilon;m)=\frac {2}{1+e^{-\frac {\pi}{\sqrt{3}}}} k^{\frac{k^2}{k^2-1}} e^{\frac{\pi}{\sqrt{3}}\frac{1-2\varepsilon}{k^ 2-1}}\sqrt{\varepsilon} +O(\varepsilon).
	\]
    \end{itemize}
  \end{itemize}
  \item [b)] Under the subcritical condition  $m=\sqrt{\varepsilon}$ we obtain that:
    \begin{itemize} 
      \item [b-1)] for $k<1$, the function $R_{3z}(h,k,\sqrt{\varepsilon};m)$ behaves as it is represented in Figure \ref{dib:FunsSyh}(c). More specifically,
      \begin{itemize}
	\item [b-1-1)] $\lim_{h\searrow h_s}R_{3z}(h,k,\varepsilon;m)>0$ and $R_{3z}(h_M,k,\sqrt{\varepsilon};m)=-e^{\frac {2\pi}{\sqrt{3}}}$, 
	\item [b-1-2)] let $h^*>0$ be a zero of $R_{3z}(h,k,\varepsilon;m)=0$, then $ \left. \frac {\partial R_{3z}}{\partial h} \right|_{(h^*,k,\sqrt{\varepsilon}:m)} <0.$
	\item [b-1-3)] denoting by $h^*(k,\varepsilon;m)$ the unique zero of $R_{3z}(h,k,\varepsilon;m)=0$ in $(h_s,h_M]$, then
	\[
	  h^*(k,\varepsilon;m)=\frac {2}{1+e^{\frac {\pi}{\sqrt{3}}}} k^{\frac{k^2}{k^2-1}} e^{\frac{\pi}{\sqrt{3}}\frac{1-2\varepsilon}{1-k^ 2}}\sqrt{\varepsilon} +O(\varepsilon),
	\]
    \end{itemize}    
      \item [b-2)] for $k\geq 1$, then $R_{3z}(h,k,\varepsilon;m)>0$ when $h\in(h_s,h_M]$.
    \end{itemize}
 \end{itemize}
\end{proposition}
\begin{proof}
By straightforward computations we write 
\begin{align*}
R_{3z}(0,k,\sqrt{\varepsilon};m)
 &= \left(1 + \frac {m+\lambda_L^s}{\lambda_L^q-\lambda_L^s}  \right)^{\frac k{\lambda_L^s}} 
  \left(1 +  \frac {m+\lambda_R^s}{\lambda_R^q-\lambda_R^s} \right)^{\frac 1{\lambda_R^s}}-
  e^{m \tau_C}\\
 &= \left(1 + \frac 1 {\frac {\lambda_L^q-\lambda_L^s}{m+\lambda_L^s}}  \right)^{\frac k{\lambda_L^s}} 
  \left(1 +  \frac 1 {\frac {\lambda_R^q+m}{-(\lambda_R^s+m)}} \right)^{-\frac 1{\lambda_R^s}}-
  e^{m \tau_C}\\
 &=\left(1+\frac 1{z_1}\right)^{\frac {k}{\lambda_L^S}}
 \left(1+\frac 1{z_2}\right)^{-\frac {1}{\lambda_R^S}},
\end{align*}
where 
\begin{align*}
 z_1&=\frac {\lambda_L^q-\lambda_L^s}{m+\lambda_L^s}=\frac{k}{m}-1+O(m),\\
 z_2&=\frac {\lambda_R^q+m}{-(m+\lambda_R^s)} = \frac{1}{m}-O(m),
\end{align*}
tend to $\infty$ (resp. to $-\infty$)  as $\varepsilon$ tends to zero when $m=\sqrt{\varepsilon}$ (resp. $m=-\sqrt{\varepsilon}$. Therefore 
\[
 \lim_{\varepsilon \searrow 0} R_{3z}(0,k,\sqrt{\varepsilon};m)= e^{\lim_{\varepsilon \searrow 0} \frac 1 {z_1} \frac k {\lambda_L^s}-\frac 1 {z_2} \frac 1 {\lambda_R^s} }-e^{\frac {2\pi}{\sqrt{3}}}=e^{\lim_{m \searrow 0} \frac{k+1}{m}+1+O(m)}-e^{\frac {2\pi}{\sqrt{3}}}.
\]
Since $h_s$ tends to zero when $\varepsilon$ does, we conclude that $R_{3z}(h_s,k,\sqrt{\varepsilon};m)<0$ when $m=-\sqrt{\varepsilon}$ and $R_{3z}(h_s,k,\sqrt{\varepsilon};m)>0$ when $m=\sqrt{\varepsilon}$, provided $\varepsilon$ is small enough. 

Expanding the different operands in the expression \eqref{cond1_3z} in power series of $\varepsilon$, and keeping the lower order terms, we obtain the following approximation of  $R_{3z}(h,k,\varepsilon;m)$ which is valid for $\varepsilon$ small enough
\begin{align*}
 R_{3z}(h,k,\varepsilon;m)
 &\approx 
 \left(\frac{h\left(1+e^{s\frac{\pi}{\sqrt{3}}}\right)}{2k\sqrt{\varepsilon}} \right)^{\frac {k^2}{\varepsilon}} 
 \left(\frac{2 e^{s\frac{\pi}{\sqrt{3}}} \sqrt{\varepsilon} }{h\left(1+e^{s\frac{\pi}{\sqrt{3}}}\right)} \right)^{\frac {1}{\varepsilon}}-
 e^{s\frac {2\pi}{\sqrt{3}}}\\
 &\approx
 \left( \frac {2\sqrt{\varepsilon}}{h \left(1+e^{s\frac{\pi}{\sqrt{3}}}\right)}  \right)^{\frac {1-k^2}{\varepsilon}}
 \left( \frac {e^{s\frac{\pi}{\sqrt{3}}}}{k^{k^2}} \right)^{\frac 1{\varepsilon}}
 -e^{s\frac {2\pi}{\sqrt{3}}},
\end{align*}
where $s$ is the $\rm{sign}(m)$. From here, we obtain the sign of the function $R_{3z}$ at $h_M=\Phi_{3z}(-1)$, depending on the parameters. That is $R_{3z}(h_M,k,\sqrt{\varepsilon};m)>0$ when $k>1$, $R_{3z}(h_M,k,\sqrt{\varepsilon};m)<0$ when $k<1$, and when $k=1$ the sign of the function $R_{3z}(h_M,k,\sqrt{\varepsilon};m)$ is equal to the sign of $m$. 

From the previous assertions, we conclude the existence of a zero $h^*\in (h_s,h_M]$ of the function $R_{3z}(h,k,\sqrt{\varepsilon};m)$ for the parameters $m=-\sqrt{\varepsilon}$ and $k<1$, and for $m=\sqrt{\varepsilon}$ and $k>1$. The expression of $h^*$  stated in the theorem follows by equalizing to zero the approximation of $R_{3z}(h,k,\sqrt{\varepsilon};m)$ given above.

The partial derivative of $R_{3z}(h,k,\varepsilon;m)$ with respect to $h$ can be written as follows 
\begin{align*}
\left.\frac {\partial R_{3z}}{\partial h}\right|_{(h,k,\sqrt{\varepsilon};m)}&=
\left( R_{3z}(h,k,\varepsilon;m)+ e^{s\frac {2\pi}{\sqrt{3}}}\right)\\
&\hspace{0.4cm}\left( 
   \frac{k}
   {\lambda_L^s((\lambda_L^q+m)(\sqrt{\varepsilon}+a)+h)}
  -\frac{1}
  {\lambda_R^s((\lambda_R^q+m)(\sqrt{\varepsilon}-a)-h)}
\right)\\
&=\left( R_{3z}(h,k,\varepsilon;m)+ e^{s\frac {2\pi}{\sqrt{3}}}\right)
\left(
  \frac{k^2-1}{h m^2}+O(\varepsilon^{-\frac 1 2})
\right).
\end{align*}
Therefore, assuming the existence of a zero $h^*$ of $R_{3z}(h,k,\sqrt{\varepsilon};m)$, we obtain that  $\left.\frac {\partial R_{3z}}{\partial h}\right|_{(h^*,k,\sqrt{\varepsilon})}<0$ if $k<1$, and $\left.\frac {\partial R_{3z}}{\partial h}\right|_{(h^*,k,\sqrt{\varepsilon})}>0$ if $k>1$, which  implies the uniqueness of such zero. From here, we conclude that function $R_{3z}$ does not change sign in $(h_s,h_M]$ neither when $m=-\sqrt{\varepsilon}$ and $k\leq 1$ nor when $m=\sqrt{\varepsilon}$ and $k\geq 1$.
\end{proof}

\subsubsection{Hyperbolicity/non-hyperbolicity of canard cycles with head}

As in the previous section, we start this one by defining the Poincaré map in neighborhood of orbits visiting the four regions $LL$, $L$, $C$ and $R$. Consider the Poincaré half-maps $\Pi_{Cd}$, $\Pi_{Cu}$, $\Pi_R$, and the time of flight $\tau_{Cd}$, $\tau_{Cu}$, and $\tau_{R}$ previously defined.

Let $\mathbf{p}_1=(-\sqrt{\varepsilon},y_1)$ be a point in the switching line $\{x=-\sqrt{\varepsilon}\}$ and located below the point $\mathbf{p}_L$, and assume that there exists a time of flight $\tau_{Ld}>0$ such that $x^L(\tau_{Ld};\mathbf{p}_1,\bm{\eta})=-1$ and $-1<x^L(s;\mathbf{p}_1,\bm{\eta})<-\sqrt{\varepsilon}$ for all $s\in(0,{\tau}_{Ld}).$ In such a case, we define the Poincar\'e half-map between the switching lines $\{x=-\sqrt{\varepsilon}\}$ and $\{x=-1\}$ at the point $y_1$ as ${\Pi}_{Ld}(y_1,{\bm{\eta}})=y^L\left(\tau_{Ld};\mathbf{p}_1,{\bm{\eta}}\right)$. Moreover, consider a point $\mathbf{p}_2=(-1,y_2)$ located below the point $\mathbf{p}_{LL}$ and assume that there exists a time of flight $\tau_{LL}>0$ such that $x^{LL}({\tau}_{LL};\mathbf{p}_1,\bm{\eta})=-1$ and $x^{LL}(s;\mathbf{p}_1,\bm{\eta})<-1$ for all $s\in(0,{\tau}_{LL}).$ We define the Poincar\'e half-map between the switching line $\{x=-1\}$ and itself at the point $y_2$ as $\Pi_{LL}(y_2,{\bm{\eta}})=y^{LL}\left(\tau_{LL};\mathbf{p}_2,{\bm{\eta}}\right)$. Finally, let $\mathbf{p}_3=(-1,y_3)$ be a point in the switching line $\{x=-1\}$, located over the point $\mathbf{p}_{LL}$, and assume that there exists $\tau_{Lu}>0$ such that $x^L(\tau_{Lu};\mathbf{p}_3,\bm{\eta})=-\sqrt{\varepsilon}$ and $-1<x^L(s;\mathbf{p}_3,\bm{\eta})<-\sqrt{\varepsilon}$ for all $s\in(0,{\tau}_{Lu}).$ We define the Poincar\'e half-map between the switching lines $\{x=-1\}$ and $\{x=-\sqrt{\varepsilon}\}$ at the point $y_3$ as ${\Pi}_{Lu}(y_3,{\bm{\eta}})=y^L\left(\tau_{Lu};\mathbf{p}_3,{\bm{\eta}}\right)$.

At this point, the Poincar\'{e} map for an orbit of system \eqref{systempal}-\eqref{vectorfield4z} 
visiting zones $LL$, $L$, $C$ and $R$ can be defined.
\begin{definition}\label{defipoincareydesplazamiento}
The Poincar\'e map $\Pi$ in the neighborhood of an orbit $\Gamma_{x_0}$ of system \eqref{systempal}-\eqref{vectorfield4z} 
visiting zones $LL$, $L$, $C$ and $R$  is defined as
\begin{equation*}
\Pi(y_0,\bm{\eta})=\Pi_{R}(\Pi_{Cu}(\Pi_{Lu}(\Pi_{LL}(\Pi_{Ld}(\Pi_{Cd}(y_0,\bm{\eta}),\bm{\eta}),\bm{\eta}),\bm{\eta}),\bm{\eta}),\bm{\eta}),
\end{equation*}
provided the composition of Poincar\'e half-maps is possible, where $y_0=\Pi_{Cd}^{-1}(\Pi_{Ld}^{-1}(\Phi_{4z}(x_0),{\bm \eta}),{\bm \eta})$.
\end{definition}

For $\varepsilon$ fixed and small enough, suppose the existence of a canard limit cycle with head $\Gamma_{x_0}$, see Figure \ref{dib:geommetry}, obtained under the parameter relation $a=\hat{a}(k,\varepsilon,x_0; m)$ given in Theorem \ref{th:exist_Gamma}. 
The cycle  $\Gamma_{x_0}$ corresponds to the fixed point of the Poincar\'e map $\Pi(y_0,{\bm  \eta}),$ where $y_0=\Pi_{Cd}^{-1}(\Pi_{Ld}^{-1}(\Phi_{4z}(x_0),{\bm \eta}),{\bm \eta})$.

The non-hyperbolicity of $\Gamma_{x_0}$ can be obtained, similarly as in the case of headless canard cycles, through the sum of the products of the traces of the matrices of the differential linear systems, and the corresponding time of flight, namely $\tau_{Cd},\tau_{Ld},\tau_{LL}$, $\tau_{Lu}$, $\tau_{Cu}$, and $\tau_{R}$.

\begin{figure}[hb]
    \centering
    \includegraphics{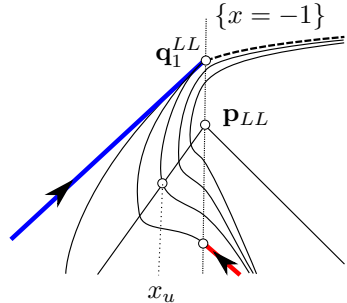}
    \begin{picture}(0,0)
    \put(-78,-7){$x_u$}
    \put(-75,83){$\mathbf{q}_1^{LL}$}
    \put(-50,60){$\mathbf{p}_{LL}$}
    \put(-53,97){$\{x=-1\}$}
    \end{picture}
    \caption{Zoom of the flow in a neighborhood of the contact point $\mathbf{p}_{LL}$. Orbits having width in $x_u<x_0<-1$ do not pass exponentially close of $\mathbf{q}_1^{LL}$. Therefore, the time of flight $\tau_{LL}$ can not be computed as in Lemma \ref{lem:times}.}
    \label{fig:Zoom2}
\end{figure}

By analogous arguments than those for headless canard limit cycles, we conclude that $\tau_{Cd}
=\tau_C(k,\varepsilon;m)$ obtained in Theorem \ref{th:connection}, and that the values of $\tau_{Cu}$ and $\tau_{Lu}$ are negligible. On the other hand, when $x_0\in (x_r,x_u)$ the canard cycle $\Gamma_{x_0}$ intersects the switching line $\{x=-1\}$ exponentially close to $\mathbf{q}_1^{LL}$, see Lemma \ref{lem:dom} and Lemma \ref{lem:Pmaps}. Therefore, the value of $\tau_{LL}$ can be approximated by the time of flight of $\Gamma_{x_0}$ from the point $(-1,h)$, where $h=\Phi_{4z}(x_0)$, to the point $\mathbf{q}_1^{LL}$. Then, from Lemma \ref{lem:times}, we obtain $\tau_{LL}=\tau_{LL}(h)$ and 
\begin{equation}\label{eq:tau_RR}
\tau_R=\tau_R(h_0)=- \frac1{\lambda_R^s} 
  \ln\left(
      1+ 
      \frac{ \lambda_R^s(\sqrt{\varepsilon}-a)+\lambda_{LL}^s(1+a)+k(\sqrt{\varepsilon}-1)+2m\sqrt{\varepsilon} }
      { (\lambda_R^q-\lambda_R^s)(\sqrt{\varepsilon}-a)}
    \right),
\end{equation}
where $h_0=-\lambda_{LL}^s(1+a)-k(\sqrt{\varepsilon}-1)-m(\sqrt{\varepsilon}+a)$ is the second coordinate of $\mathbf{q}_1^{LL}$.
Notice that when $x_u<x_0<-1$, we can not assure that $\Gamma_{x_0}$ intersects $\{x=-1\}$ exponentially close to $\mathbf{q}_1^{LL}$, see Figure \ref{fig:Zoom2}.
In this case, neither expression in Lemma \ref{lem:times} nor expression \eqref{eq:tau_RR} are good approximations for $\tau_{LL}$ and $\tau_R$, respectively. Therefore, we have eliminated the interval $(x_u,-1)$ from the stated of the Theorem \ref{th:supercritico} and of the Theorem \ref{th:subcritico}.

Following similar arguments that those applied in Section \ref{sec:3z}, a necessary condition on the canard cycle $\Gamma_{x_0}$ to be non-hyperbolic can be written as $R_{4z}(h,k,\varepsilon;m)=e^{t_L \tau_{Ld} + t_{LL} \tau_{LL} + t_R \tau_{R}}-e^{m \tau_C }=0$. By using the expressions of $\tau_{Ld}, \tau_{LL},\tau_{R}$ and $\tau_C$ we obtain  
\begin{align}\label{cond1_4z}
R_{4z}(h,k,\varepsilon;m)=  
   &\left(
      1 + 
      \frac {h+m(\sqrt{\varepsilon}+a)+\lambda_L^s(2\sqrt{\varepsilon}+a-1)}
      {(\lambda_L^q-\lambda_L^s)(\sqrt{\varepsilon}+a)}  
    \right)^{\frac k{\lambda_L^s}}\nonumber\\
   &\left(
      1+ 
      \frac{ h+m(\sqrt{\varepsilon}+a)+k(\sqrt{\varepsilon}-1)+\lambda_{LL}^s(1+a)}
           { (\lambda_{LL}^q-\lambda_{LL}^s)(1+a) }
    \right)^{\frac1{\lambda_{LL}^s}}\\
   &\left(
      1+ 
      \frac{ \lambda_R^s(\sqrt{\varepsilon}-a)+\lambda_{LL}^s(1+a)+k(\sqrt{\varepsilon}-1)+2m\sqrt{\varepsilon} }
      { (\lambda_R^q-\lambda_R^s)(\sqrt{\varepsilon}-a)}
    \right)^{\frac1{\lambda_R^s}}-e^{m\tau_C}.\nonumber  
\end{align}
\begin{proposition}\label{prop:R4z}
For $\varepsilon$ fixed and small enough, there exists $0<\delta\ll 1,$ such that, for  $x_0\in(x_r,x_u)$ and $h=\Phi_{4z}(x_0)$:
\begin{itemize}
 \item [a)] if $R_{4z}(h,k,\sqrt{\varepsilon}; m)<-\delta$, then $\Gamma_{x_0}$ is a hyperbolic stable canard cycle with head;
 \item [b)] if $h$ is a simple root of $R_{4z}(h,k,\sqrt{\varepsilon}; m)$, then in a neighborhood of $\Gamma_{x_0}$ there is a nonhyperbolic canard cycle with head;
 \item [c)] if $R_{4z}(h,k,\sqrt{\varepsilon}; m)>\delta$, then $\Gamma_{x_0}$ is a hyperbolic unstable canard cycle with head.
\end{itemize}
\end{proposition}
\begin{proof}
The proposition follows similarly to Proposition \ref{prop:R3z}.
\end{proof}

Next, we describe the qualitative behavior of  $R_{4z}(h,k,\varepsilon;m)$, as a function of $h$, for fixed values of the parameters $k$ and $\varepsilon$.  Even when the domain of the function $R_{4z}$ is greater, we consider it reduced to $(h_r,h_u)$, where $h_r=\Phi_{4z}(x_r)$ and $h_u=\Phi_{4z}(h_u)$, see Lemma \ref{lem:dom}.

\begin{figure}[ht]
\begin{center}
\includegraphics{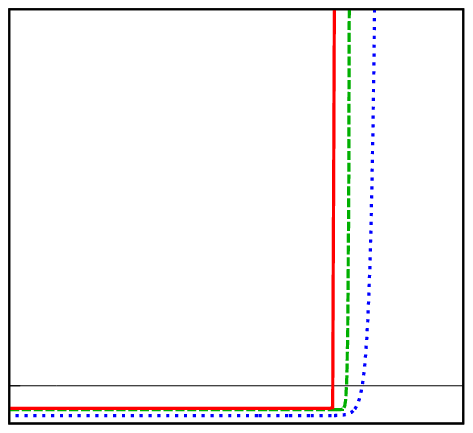}\quad\quad\includegraphics{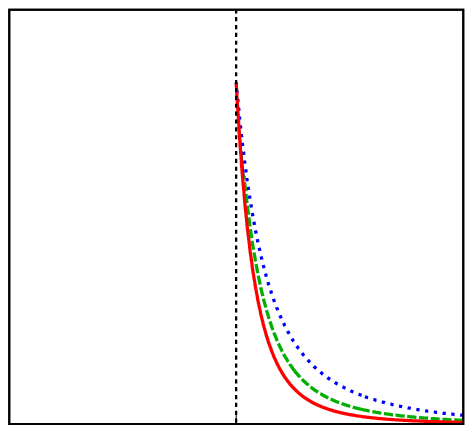}\\ 
\includegraphics{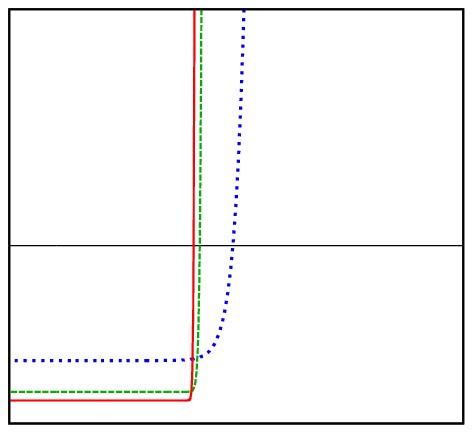}\quad\quad\includegraphics{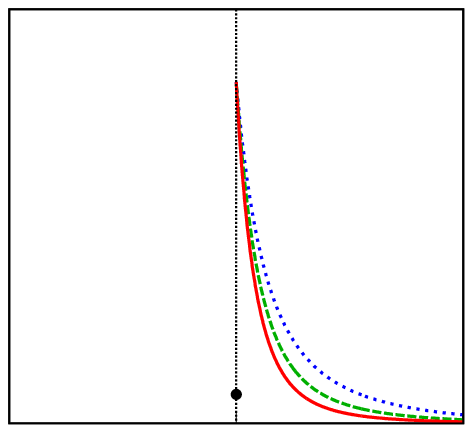}
\begin{picture}(0,0)
\put(-360,275){(a)}
\put(-236,254){\color{white}\circle*{15}}
\put(-246,254){\color{white}\circle*{15}}
\put(-280,250){{\small $R_{4z}(h,k,\sqrt{\varepsilon};m)$}}
\put(-220,147){$h$}
\put(-352,167){$0$}
\put(-246,171){\circle*{4}}
\put(-242,171){\circle*{4}}
\put(-237,171){\circle*{4}}
\put(-259,175){$h^*$}
\put(-170,275){(b)} 
\put(-86,147){$1$}
\put(-30,147){$k$}
\put(-165,245){$h^*$}
\put(-60,147){$1.3$}
\put(-55,160){\line(0,1){25}}
\put(-55,162){\circle*{4}}
\put(-55,166){\circle*{4}}
\put(-55,171){\circle*{4}}
\put(-61,204){\color{white}\circle*{15}}
\put(-73,215){{\small $h^*(k,\sqrt{\varepsilon};m)$}}
\put(-360,130){(c)}
\put(-352,65){$0$}
\put(-274,85){\color{white}\circle*{15}}
\put(-280,82){{\small $R_{4z}(h,k,\sqrt{\varepsilon};m)$}}
\put(-220,4){$h$}
\put(-170,130){(d)}
\put(-86,4){$1$}
\put(-30,4){$k$}
\put(-165,100){$h^*$}
\put(-104,64){\color{white}\circle*{15}}
\put(-145,61){$h^*(k,\sqrt{\varepsilon};m)$}
\end{picture}\end{center}
\caption{Graphical representation of the functions $R_{4z}(h,k,\sqrt{\varepsilon};m)$ and $h^*(k,\varepsilon;m)$ analyzed in Theorem \ref{th:nhcond4z} for different values of $\varepsilon$. In particular blue/dotted curves correspond with $\varepsilon=0.05$, green/dashed curves with  $\varepsilon=0.01$ and red/solid curves with $\varepsilon=0.005$. First column contains the graphs of the function $R_{4z}(h,k,\varepsilon;m)$ as a function of $h$ for the previous values of $\varepsilon$: panel (a) when $k=1.3$ and $m=-\sqrt{\varepsilon}$ and panel (c) when $k=0.75$ and $m=\sqrt{\varepsilon}$. Second column contains the graphs of the function $h^*(k,\varepsilon;m)$ as a function of $k$ for the previous values of $\varepsilon$: panel (b) when $m=-\sqrt{\varepsilon}$ and panel (c) when $m=\sqrt{\varepsilon}$.}\label{fig:FunS4z}
\end{figure}

\begin{proposition}\label{th:nhcond4z}
Fixed $\varepsilon$ small enough, we consider the function $R_{4z}(h,k,\varepsilon;m)$ defined in \eqref{cond1_4z}. 
\begin{itemize}
      \item [a)] If $k<1$, or $k=1$ and $m=-\sqrt{\varepsilon}$, then $R_{4z}(h,k,\sqrt{\varepsilon}; m)<0$ in $(h_r,h_u)$.
      \item [b)] If $k>1$, or $k=1$ and $m=\sqrt{\varepsilon}$, then $R_{4z}(h,k,\sqrt{\varepsilon};m)$ behaves as in Figure \ref{fig:FunS4z}(a) or in Figure \ref{fig:FunS4z}(c), depending on the supercritical case, $m=-\sqrt{\varepsilon}$, or the subcritical case, $m=\sqrt{\varepsilon}$, respectively. More specifically: 
	\begin{itemize}
	  \item [b-1)] $R_{4z}(h_r,k,\sqrt{\varepsilon};m)<0$ and $R_{4z}(h_u,k,\sqrt{\varepsilon};m)>0$, see Figure \ref{fig:FunS4z}(a) and (c);
	  \item [b-2)] let $h^*\in(h_r,h_u)$ be a zero of $R_{4z}(h,k,\varepsilon;m)$, then $ \left. \frac {\partial R_{4z}}{\partial h} \right|_{(h^*,k,\sqrt{\varepsilon};m)} >0$; 
	  \item [b-3)] denoting by $h^*(k,\sqrt{\varepsilon};m)$ the unique zero of $R_{4z}(h,k,\varepsilon;m)$ in $(h_r,h_u)$, it follows that, if $k>1$, then
	      \[
		  h^*(k,\sqrt{\varepsilon};m)=(k+1)e^{\frac {2-k^2}{2}} (\sqrt{\varepsilon})^{\frac {k^2-1}{k^2}}+O(\sqrt{\varepsilon}),
	      \]
	      see Figure \ref{fig:FunS4z}(b) and (d), and if $k=1$ then 
	      \[
		  h^*(1,\sqrt{\varepsilon};m)= \frac {2}{1+e^{\frac{\pi}{\sqrt{3}}}}+O(\sqrt{\varepsilon}),
	      \]	      
	      see Figure \ref{fig:FunS4z}(d).
	\end{itemize}
 \end{itemize}
\end{proposition}
\begin{proof}
  Setting $m=\pm \sqrt{\varepsilon}$ and $s=\rm{sign}(m)$, and expanding in power series of $\varepsilon$ every term in the expression of $R_{4z}(h,k,\sqrt{\varepsilon}; m)$, we obtain that 
 \begin{align*}
  R_{4z}(h,k,\sqrt{\varepsilon}; m) =& 
    \left(
      1+k-h
      -\frac{2k+h(e^{s\frac{\pi}{\sqrt{3}}}-1)}{1+e^{s \frac{\pi}{\sqrt{3}}}}
      \sqrt{\varepsilon}
      +O(\varepsilon)
    \right)^{-\frac{1}{\varepsilon}+1+O(\varepsilon)}\\
   &\left(
      \frac{k(1+e^{-s\frac{\pi}{\sqrt{3}}})}{2}
      \frac 1{\sqrt{\varepsilon}}
      -\frac{k(1+e^{-s\frac{\pi}{\sqrt{3}}})-2}{2}
      +O(\varepsilon^{\frac 1 2})
    \right)^{-\frac{1}{\varepsilon}+1+O(\varepsilon)}\\
   &\left(
      \frac{h(1+e^{s\frac{\pi}{\sqrt{3}}})}{2k}
      \frac 1{\sqrt{\varepsilon}}
      +1
      -\left(
	\frac{e^{s\frac{\pi}{\sqrt{3}}}+1-2ks}{2k^2}
	+\frac{e^{s\frac{\pi}{\sqrt{3}}}(k^2-5)-4}{4k^3}h
      \right)\sqrt{\varepsilon}
      +O(\varepsilon)
    \right)^{\frac{k^2}{\varepsilon}-1+O(\varepsilon)}\\
   &-e^{s\frac{2\pi}{\sqrt{3}}-s\frac{k+1}{k}\sqrt{\varepsilon}+O(\varepsilon)}.
 \end{align*}
By changing the sign in the exponent of the first two terms, it follows that
 \begin{align}\label{ex:TayS}
 R_{4z}(h,k,\sqrt{\varepsilon}; m) =&
    \left(
      \frac{1}{1+k-h}
     +\frac{ 2k+h(e^{s\frac{\pi }{\sqrt{3}}}-1) }{\left( 1+ e^{s\frac{\pi }{\sqrt{3}}} \right) {{\left(1+k-h\right) }^{2}}}\sqrt{\varepsilon} +O(\varepsilon)
    \right)^{\frac{1}{\varepsilon}-1+O(\varepsilon)}\nonumber \\ 
   &\left(
      \frac{2}{k(1+e^{-s\frac{\pi}{\sqrt{3}}})}
      {\sqrt{\varepsilon}}
      +O(\varepsilon)
    \right)^{\frac{1}{\varepsilon}-1+O(\varepsilon)}\\ 
      &\left(
      \frac{h(1+e^{s\frac{\pi}{\sqrt{3}}})}{2k}
      \frac 1{\sqrt{\varepsilon}}
      +1
      -\left(
	\frac{e^{s\frac{\pi}{\sqrt{3}}}+1-2ks}{2k^2}
	+\frac{e^{s\frac{\pi}{\sqrt{3}}}(k^2-5)-4}{4k^3}h
      \right)\sqrt{\varepsilon}
      +O(\varepsilon)
    \right)^{\frac{k^2}{\varepsilon}-1+O(\varepsilon)} \nonumber \\
   &-e^{s\frac{2\pi}{\sqrt{3}}-s\frac{k+1}{k}\sqrt{\varepsilon}+O(\varepsilon)}. \nonumber
 \end{align}
 Then
 \begin{align*}
  R_{4z}(0,k,\sqrt{\varepsilon};m) =& 
    \left(
      \frac{1}{1+k} +O(\sqrt{\varepsilon})
    \right)^{\frac{1}{\varepsilon}-1+O(\varepsilon)}
   \left(
      \frac{2}{k(1+e^{-s\frac{\pi}{\sqrt{3}}})}
      \sqrt{\varepsilon}
      +O(\varepsilon)
    \right)^{\frac{1}{\varepsilon}-1+O(\varepsilon)}\\
   &\left(
      1
      -\frac{e^{s\frac{\pi}{\sqrt{3}}}+1-2sk}{2k^2}
      \sqrt{\varepsilon}
      +O(\varepsilon)
    \right)^{\frac{k^2}{\varepsilon}-1+O(\varepsilon)}\\
   &-e^{s\frac{2\pi}{\sqrt{3}}-s\frac{k+1}{k}\sqrt{\varepsilon}+O(\varepsilon)},
 \end{align*}
which can be approximated, for $\varepsilon$ small enough, as follows 
\[
  R_{4z}(0,k,\sqrt{\varepsilon};m) \approx
    \left(
      \frac{2\sqrt{\varepsilon}}{k(1+k)(1+e^{-s\frac{\pi}{\sqrt{3}}})}
      -\frac{(e^{s\frac{\pi}{\sqrt{3}}}+1-2sk)\varepsilon}{k(1+k)(1+e^{-s\frac{\pi}{\sqrt{3}}})}
      + O(\varepsilon^{\frac 3 2})
    \right)^{\frac 1 {\varepsilon}}
   -e^{s\frac{2\pi}{\sqrt{3}}-s\frac{k+1}{k}\sqrt{\varepsilon}+O(\varepsilon)}.
\]
First operand in the right side of previous equation tends to zero provided $\varepsilon$ does. Hence, for $\varepsilon$ small enough we obtain that  $-e^{s\frac{2\pi}{\sqrt{3}}}< R_{4z}(0,k,\sqrt{\varepsilon}; m)<0$. Same expression is also satisfied by  $R_{4z}(h_r,k,\sqrt{\varepsilon}; m)$  since $\lim_{\varepsilon \searrow 0} h_r=0$. Thus, we conclude that $-e^{s\frac{2\pi}{\sqrt{3}}}< R_{4z}(h_r,k,\sqrt{\varepsilon}; m)<0$.  

On the other hand, from expression \eqref{ex:TayS} it follows that 
 \begin{align*}
  R_{4z}(k,k,\sqrt{\varepsilon}; m) 
  &=\left(
      1+k\sqrt{\varepsilon}+O(\varepsilon)
    \right)^{\frac{1}{\varepsilon}-1+O(\varepsilon)}
    \left(
      \frac{2}{k(1+e^{-s\frac{\pi}{\sqrt{3}}})}
      \sqrt{\varepsilon}
      +O(\varepsilon)
    \right)^{\frac{1}{\varepsilon}-1+O(\varepsilon)}\\
  &\phantom{=}\left(
      \frac{1 + e^{s\frac{\pi}{\sqrt{3}}}}{2}
      \frac 1{\sqrt{\varepsilon}} 
      +1
      +O(\sqrt{\varepsilon})
    \right)^{\frac{k^2}{\varepsilon}-1+O(\varepsilon)}-e^{s\frac{2\pi}{\sqrt{3}}-s\frac{k+1}{k}\sqrt{\varepsilon}+O(\varepsilon)}\\
   &=\left(
    \frac {2^{1-k^2}(1+e^{s\frac{\pi}{\sqrt{3}}})^{k^2-1} e^{s\frac{\pi}{\sqrt{3}}}}{k} \varepsilon^{\frac{1-k^2}{2}} + \ldots
  \right)^{\frac 1 {\varepsilon}}-e^{s\frac{2\pi}{\sqrt{3}}}.
 \end{align*}
Hence, for $\varepsilon$ small enough we obtain that if $k>1$, then $R_{4z}(k,k,\sqrt{\varepsilon}; m)>0$; and if $k<1$, then $R_{4z}(k,k,\sqrt{\varepsilon}; m)<0$. Moreover, when $k=1$ it follows that $R_{4z}(k,k,\sqrt{\varepsilon}; m)$ has the same sign than $m$. Since $\lim_{\varepsilon \searrow 0} h_u=k$, previous inequalities are also satisfied by  $R_{4z}(h_u,k,\sqrt{\varepsilon}; m)$.
 
Let $h_r<h^*<h_u$ be a zero of $R_{4z}(h,k,\sqrt{\varepsilon}; m)$. Straight forward computations shows that
 \begin{align*}
  \left.\frac{\partial R_{4z}}{\partial h}\right|_{(h^*,k,\sqrt{\varepsilon};m)}&=
     e^{s\frac{2\pi}{\sqrt{3}}-s\frac{k+1}{k}\sqrt{\varepsilon}+O(\varepsilon)}
     \left(
	\frac{ k^3+k^2(1-h^*)+h^*}
	     {h^*(k+1-h^*)\varepsilon} 
	 +O\left(\frac 1 {\sqrt{\varepsilon}}\right)
     \right)>0,
 \end{align*}
which implies that, when exists, the zero $h^*$ is unique. Moreover, by keeping the lower order terms in \eqref{ex:TayS}, we obtain the following implicit expression for an approximation to the solution $h^*(k,\sqrt{\varepsilon})$ of $R_{4z}(h,k,\sqrt{\varepsilon}; m)=0$, that is
 \[
  \left( 
    \frac 1{1+k-h^*} 
  \right)
  \left( 
    \frac {2\sqrt{\varepsilon}}{k(1+e^{-s\frac {\pi}{\sqrt{3}}})}
  \right)
  \left(
    \frac {h^*(1+e^{s\frac {\pi}{\sqrt{3}}})}{2k} \frac 1{\sqrt{\varepsilon}}  
  \right)^{k^ 2}
  =e^{s \frac {2\pi}{\sqrt{3}} \varepsilon},
 \]
When $k=1$, the equation above can be solved and it follows that 
\[
 h^*= \frac {2}{1+e^{s\frac{\pi}{\sqrt{3}}}}+O(\sqrt{\varepsilon}).
\]
which implies that $h^*<h_u$ only if $m=\sqrt{\varepsilon}$.
When $k\neq 1$ the solution of the equation can be approximated by the undetermined coefficients method
\[
 h^*=(k+1)e^{\frac {2-k^2}{2}} (\sqrt{\varepsilon})^{\frac {k^2-1}{k^2}}+O(\sqrt{\varepsilon}).
\]
\end{proof}

In Figure \ref{fig:FunS4z}(a)  we represent $R_{4z}(h,k,\varepsilon;m)$ as a function of $h$ for different values of $\varepsilon$ and $k>1$. We notice that function $R_{4z}(h,k,\varepsilon;m)$ is negative for $h$ close to $h_r$, and positive for values close to $h_{u}$. Moreover, as it follows from the proof of the previous theorem, $\lim_{\varepsilon \searrow 0} R_{4z}(0,k,\sqrt{\varepsilon};m)=-e^{\frac {2\pi}{\sqrt{3}}}$. Furthermore, the unique zero $h^*(k,\sqrt{\varepsilon})$ tends to zero as $\varepsilon$ tends to zero. In Figure \ref{fig:FunS4z}(b) we represent function $h^*(k,\sqrt{\varepsilon})$ as a function of $k$ for different values of $\varepsilon$. Notice that all the curves have a common point at $k=1$ and $h^*=2{e}^{\frac 1 2}\approx 3.297442541400256\ldots>k,$ which lies out of the interval $(h_r,h_u)$ since $h_u\approx k$. Therefore this part of the curves do not correspond to non-hyperbolic canard cycles, see Figure \ref{dib:hk34}.


\subsubsection{Correspondence to saddle-node bifurcations}
To finish with the proof of the theorems \ref{th:supercritico} and \ref{th:subcritico}, in this subsection we are going to prove that the non-hyperbolic canard limit cycles whose existence has been proved in the previous sections correspond, indeed, to a saddle-node bifurcation. We focus the proof in the case of headless canards. The proof in the case of canards with head is analogous.

We will prove that the non-degeneracy conditions on the Poincar\'e map hold, that is,  the second derivative of the Poincar\'e map   with respect to the initial condition and the derivative of the Poincar\'e map with respect to parameter $a$ are both different from zero, see \cite{K04,Perko}. 

Consider a non-hyperbolic headless canard cycle $\Gamma_{x_0}$ corresponding to a fixed point of the Poincar\'e map $y_0$
for parameters $\bm{\eta}=\bar{\bm{\eta}},$ i.e. 
$y_0=\Pi_{Cd}^{-1}(\Pi_{L}^{-1}(\Phi(x_0),\bar{{\bm \eta}}),\bar{{\bm \eta}})$.
First,  we compute the second derivative of the Poincar\'e map with respect to the initial condition and we will see that, in a headless non-hyperbolic canard cycle this derivative is nonzero if and only function $R_{3z}$ is non zero at $\bar{h}=\Phi_{3z}(x_0)$. Second, we compute the derivative of the Poincar\'e map with respect to parameter $a$ and we test that this derivative does not vanish. 


 From 
expression (\ref{derprim}),
the second derivative of Poincar\'e map with respect to $h$ in the non-hyperbolic fixed point $\bar{h}$ is given by
\begin{equation}\label{derseg}
\frac{\partial^2 \Pi}{\partial y^2}(y_0,\bar{\bm{\eta}})=\frac{\partial G}{\partial h}(\bar{h},\bar{\bm{\eta}}), 
\end{equation}
where $$G({h},\bm{\eta})=e^{t_L \tau_L-m\tau_C+t_R\tau_R}.$$

It is easy to see that 
$$G({h},\bm{\eta})=(R_{3z}(h,k,\varepsilon;m)+e^{m\tau_C})e^{-m\tau_C}=R_{3z}(h,k,\varepsilon;m)e^{-m\tau_C}+1.$$
As in a non-hyperbolic headless canard cycle $R_{3z}(h,k,\varepsilon;m)=0,$ it is easy to see that,
$$\frac{\partial G}{\partial h} (\bar{h},\bar{\bm{\eta}})=\frac{\partial R_{3z}}{\partial h} (\bar{h},\bar{\bm{\eta}})e^{-m\tau_C}.$$
Thus, from expression (\ref{derseg}), in a non-hyperbolic headless canard cycle,
\begin{equation}\label{dersegfin}
\frac{\partial^2 \Pi}{\partial y^2}(y_0,\bar{\bm{\eta}})=\frac{\partial R_{3z}}{\partial h} (\bar{h},\bar{\bm{\eta}})e^{-m\tau_C}, 
\end{equation} 
and then the sign of $\frac{\partial^2 \Pi}{\partial y^2}(y_0,\bar{\bm{\eta}})$ is that of $\frac{\partial R_{3z}}{\partial h} (\bar{h},\bar{\bm{\eta}}),$ which is non zero from Theorem \ref{th:nhcond3z} (a-2-2), (b-1-2) and Theorem \ref{th:nhcond4z} (b-2).

To compute the derivative of the Poincar\'e map with respect to parameter $a$ in the neighborhood of a non-hyperbolic headless canard cycle $\Gamma_{x_0}$, we will use the following reasoning. Taking a point $(\sqrt{\varepsilon},y_1)$  in a neighborhood of $(\sqrt{\varepsilon},y_0)$ where  $y_0=\Pi_{Cd}^{-1}(\Pi_{L}^{-1}(\Phi_{3z}(x_0),{\bm \eta}),{\bm \eta})$. The image through the Poincar\'e map of $y_1$, that is $\Pi(y_1,{\bm \eta})$, will be exponentially close to the second component of $ \mathbf{q}_1^R,$ that is,
\begin{equation*}\label{derpda1}
\Pi(y_1,\bm{\eta})=(m+\lambda_R^s)(\sqrt{\varepsilon}-a)+\chi(y_1,\bm{\eta}),
\end{equation*}
with function $\chi(y_1,\bm{\eta})$ and its derivatives $O(\exp(-c/\varepsilon))$ small, where $c$ is a positive constant depending on $y_1.$
Hence,
\begin{equation*}\label{derpda2}
\frac{\partial \Pi}{\partial a}(y_1,\bm{\eta})=-(m+\lambda_R^s)+\frac{\partial \chi}{\partial a}(y_1,\bm{\eta})\simeq -m+\varepsilon, 
\end{equation*}
and $\frac{\partial \Pi}{\partial a}(y_1,\bm{\eta})\neq 0$ as $m=\pm \sqrt{\varepsilon}\neq 0.$

\subsection{Proof of Theorem \ref{th:supercritico_LP}} 

In Figure \ref{dib:FunsSyh}(a) we draw the graph of the function  $R_{3z}(h,k,\sqrt{\varepsilon};m)$ as a function of $h$ by fixing parameter $k=2.5$ and parameter $\varepsilon\in\{0.05,0.01,0.005\}$. As it has been obtained in the proof of Proposition \ref{th:nhcond3z}, $R_{3z}(h,k,\sqrt{\varepsilon};m)$ tends to $ -e^{\frac {2\pi}{\sqrt{3}}}$ as $h$ decrease to $h_s$, and it tends to $\infty$ as $h$ tends to $h_M$. In panel (a) we also represent the zero, $h^*\in (h_s,h_M]$, at which the function change sign. As it can be observed this zero tends to zero as $\varepsilon$ tends to zero.  


From Theorem \ref{th:nhcond3z}(c), in Figure \ref{dib:FunsSyh}(b) we draw function $h^*(k,\sqrt{\varepsilon})$ for different values of $k<1$. As it can be observed, for a fixed $k<1$, as $\varepsilon$ tends to zero, so tends the zero $h^*(k,\sqrt{\varepsilon})$, that is  the height of the saddle-node tends to zero assuming $\varepsilon$ does. We conclude that fixed $k$ the singular limit of every saddle-node is the equilibrium point at the fold. 

Nevertheless, for any $h_0\in (h_s,h_M] \cup (h_r,h_u)$, we prove that there exists a suitable election of the parameter $k$ such that 
\[
 \lim_{\varepsilon\searrow 0} h^*(k(\varepsilon),\varepsilon;m)=h_0,
\]
which means that we can chose parameter $k$ in such a way that the singular limit of the saddle-note to be a singular cycle of prefixed height, $h_0$. 

Next, we address the supercritical case $m=-\sqrt{\varepsilon}$ and $h_0\in(h_s,h_M]$. The remainder cases follow by similar arguments. In order to obtain the parameter value $k$ allowing the existence of a saddle-node canard cycle with height $h_0$, we use implicit equation $h^*(k,\varepsilon;m)=h_0$, where the function $h^*(k,\varepsilon;m)$ is given in Proposition \ref{th:nhcond3z}(a-2-3). Since the partial derivative 
$$\frac{\partial h^*}{\partial k}(k,\varepsilon;m)=\frac{2 \sqrt{\varepsilon} k^{\frac{k^2}{k^2-1}+1}
   e^{\frac{\pi  \left(k^2-2 \varepsilon\right)}{\sqrt{3}
   \left(k^2-1\right)}} \left(2 \sqrt{3} \pi  (2
   \varepsilon-1)+3 k^2-6 \ln(k)-3\right)}{3
   \left(1+e^{\frac{\pi }{\sqrt{3}}}\right)
   \left(k^2-1\right)^2},$$
is different from zero when $k$ is greater but close to $1$ and $\varepsilon>0$, the Implicit Function Theorem can be applied, and there exists $k(h,\varepsilon;m)$ 
such that the differential system \eqref{systempal}-\eqref{vectorfield4z} with parameters $k=k(h,\varepsilon;m)$ and $a=\hat{a}(k(h,\varepsilon;m),\varepsilon,x_0;m)$ exhibits the saddle-node canard $\Gamma_{x_{0}}$, with $x_0=\Phi_{3z}^{-1}(h)$, whose existence has been stated in Theorem \ref{th:supercritico}.

\section{Conclusions and perspectives}\label{sec_conclusions}

\begin{figure}[ht]
 \begin{center}
  \includegraphics{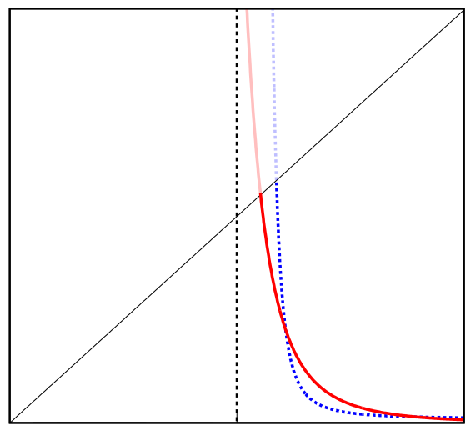}\quad\quad\includegraphics{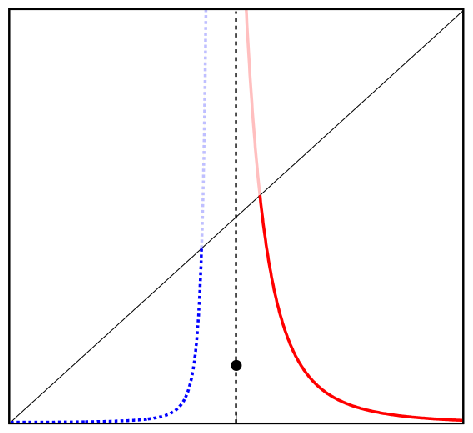}
  \begin{picture}(0,0)
   \put(-365,125){(a)}
   \put(-362,100){$h^*$}
   \put(-283,5){$1$}
   \put(-220,5){$k$}
   \put(-170,125){(b)}
   \put(-167,100){$h^*$}
   \put(-87,5){$1$}
   \put(-25,5){$k$}
   \put(-325,620){$h_M(k)$}
   \put(-140,53){$h_M(k)$}
   \multiput(-94,15)(0,5){10}{\line(0,1){2}}
   \put(-103,5){$k_1$}
   \multiput(-78,15)(0,5){14}{\line(0,1){2}}
   \put(-73,5){$k_2$}
  \end{picture}
  \end{center}
 \caption{Curves of saddle-node canard cycles of system \eqref{systempal}-\eqref{vectorfield4z} with $\varepsilon=1e-5$. Pointed/blue curves correspond with three zonal saddle-node canard cycles (headless canards). Solid/red curves correspond with four zonal saddle-node canard cycles (canards with head). The diagonal is the value of the height $h_M$ as a function of $k$. The part of the curves over the diagonal corresponds to zeros of the function $R_{3z}$ and $R_{4z}$ which are not saddle-node canards. 
 Panel (a) represents the two saddle-node canard cycles appearing in the supercritical case $m=-\sqrt{\varepsilon}$ for $k>1 $. Panel (b) represents the saddle-node canard cycles appearing in the subcritical case $m=\sqrt{\varepsilon}$, in this case for each value of $k$ only one saddle-node limit cycle can appears.}\label{dib:hk34} 
 \end{figure}

\begin{figure}
 \begin{center}
    \includegraphics[width=0.45\textwidth]{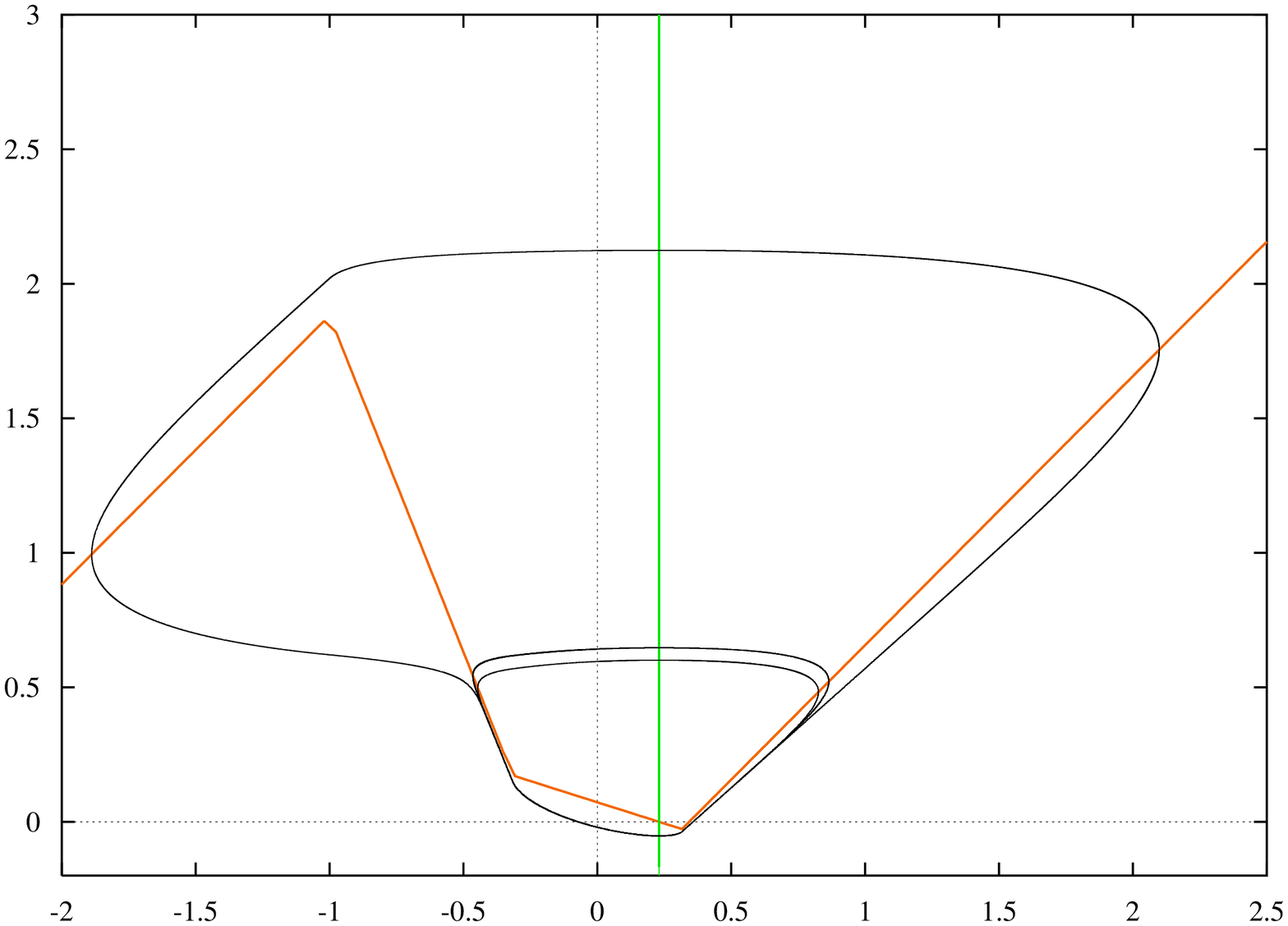}\quad\includegraphics[width=0.45\textwidth]{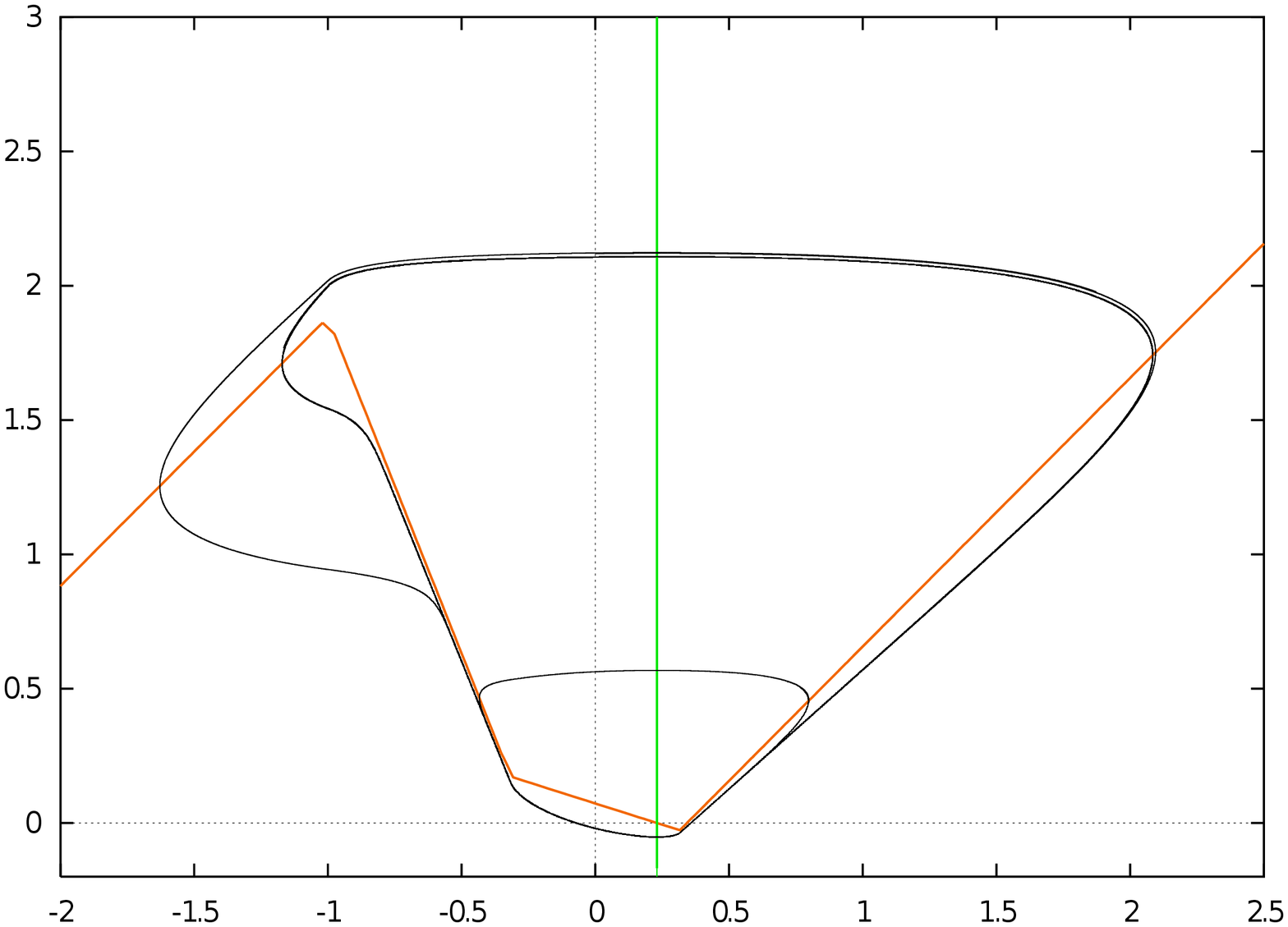}
 \end{center}
 \caption{Simulation of three canard limit cycles close to a saddle-node configuration for the supercritical case $m=-\sqrt{\varepsilon}$. According with Figure \ref{dib:hk34}(a) two configuration are possible when $k>1$: a saddle-node canard cycle without head and saddle-node canard cycle with head. Panel (a) depicts a configuration close to the saddle-node canard cycle without head for $\varepsilon=0.1$, $k=2.5$ and $a=a_1=0.2305968812$.  The three canard cycles have initial conditions $(0,0.595), (0,0.642), (0,2.12361)$. Panel (b) depicts a configuration close to the saddle-node canard cycle with head for $\varepsilon=0.1$, $k=2.5$ and $a=a_2=0.23059688315966$, the three canard cycles have initial conditions $(0,0.561), (0,2.10673), (0,2.12)$. Note that $|a_1-a_2|\approx 2e-9\approx e^{-\frac 2{\varepsilon}}$, then both saddle-node canard cycles take place exponentially close.}\label{fig3canards}
\end{figure}

In this work we have analyzed the existence of saddle-node bifurcation of canard cycles in PWL systems. We have revised the results in the smooth framework obtained in \cite{KS01}, but in the PWL context. As we have already commented, there, the authors considered two different scenarios, depending whether the Hopf bifurcation where the family of cycles is born is supercritical or subcritical. Here, we point out the similarities and differences that we have found in this study in both cases:
\begin{itemize}
\item In \cite{KS01}, canard cycles develop along a branch born at a Hopf bifurcation, at $a=a_H$, and the canard explosion takes place at a value which is at a distance of $O(\varepsilon)$ from the $a_H.$ In the PWL context, 
we have checked that the canard explosion takes place at a value which is at a distance of $O(\sqrt{\varepsilon})$ from the $a_H.$
\item \textit{Supercritical case, $m=-\sqrt{\varepsilon}$}: System \eqref{systempal}-\eqref{vectorfield4z} is able to reproduce the dynamics in the smooth case with $k \leq 1$, that is, the existence of a family of stable canard cycles.
By letting $k$ to increase, we find new scenarios that have not been reported in the smooth framework. In particular, when $k>1,$ we find situations where two saddle-node bifurcations of canard cycles take place, one of headless canards and another one of canards with head, see Figure \ref{dib:hk34}, panel (a). In this case, three canard limit cycles can coexist, see Figure \ref{fig3canards}.
\item \textit{Subcritical case, $m=\sqrt{\varepsilon}$}: In this case, system \eqref{systempal}-\eqref{vectorfield4z}
is able to reproduce the dynamics in the smooth case, with the advantage that in the PWL case we can easier control the different behaviors that appear. In particular, we have proved the existence of saddle-node bifurcation of headless canards for $k<1,$ and of canards with head for $k>1,$ see Figure \ref{dib:hk34}, panel (b).
\end{itemize}
Moreover, in Theorem \ref{th:supercritico_LP} we have stated that, both in subcritical and supercritical cases,  for every height between the smallest canard cycle and the relaxation oscillation cycle there exist parameters $k$ and $\varepsilon$ such that a saddle-node canard limit cycle with this height exists.

Furthermore, note that, in Figure \ref{dib:hk34} we have represented the height of the saddle-node canard versus the parameter $k$ for a particular value of $\varepsilon$. The straight line corresponds with the height, $h_M$, of the canard orbit through the tangent point $\mathbf{p}_{LL}$, and coincides with the maximum height of a 3-zones saddle-node canard orbit. Hence, for $k_1<k<k_2$ the saddle-node canard orbits predicted by function $h^*(k,\sqrt{\varepsilon};m)$ do not correspond neither with 3-zones nor 4-zones saddle-node canard orbits. It is clear that there must be a transition between both types of saddle-node canard orbits and, as we have mentioned along the manuscript, such a transition can not be obtained through the functions $R_{3z}$ and $R_{4z}$, see \eqref{cond1_3z} and \eqref{cond1_4z} respectively. In fact, similar situation happens in the smooth context, see \cite{DDR14}, where slow divergence integral results to be not enough for the analysis of these transitory orbits. Although the results in this article do not apply to transitory canard cycles, from them, we can establish the following conjecture:

 \begin{conjecture} For $\varepsilon=0$, consider the transitory singular canard $\Gamma$, see  \cite{DDR14}, formed by the critical manifold in $[-1,1]$ and the segment $\{(x,k):x\in[-1,1]\}$.  There exist values of $\varepsilon>0$ and $k \approx 1$ such that corresponding system \eqref{systempal}-\eqref{vectorfield4z} exhibits one, two or three canard limit cycles in a neighborhood of $\Gamma$, even more, the system also exhibits zero, one or two saddle-node canard cycles in a neighborhood of $\Gamma$.
\end{conjecture}

Finally, in the subcritical case, the existence of saddle-node bifurcation of headless canards for $k<1,$ and of canards with head for $k>1,$ leads us to come up with:
\begin{conjecture}
 System \eqref{systempal}-\eqref{vectorfield4z} in the subcritical case ($m=\sqrt{\varepsilon}$) possesses a saddle-node bifurcation of transitory canard cycles  for $k=1.$ 
\end{conjecture}

The use of this simpler family of slow-fast systems to reproduce canard dynamics bring us some information which could be interesting when revisiting the smooth context. In particular, the conditions $k<1$ and $k>1$ organizing  the dynamics in the main results, suggest the importance of the ratio between the slopes of the fast nullcline in order to exhibit or not  saddle-node canard cycles with head. Assuming this idea, it seems that only saddle-node canard cycles with head can appear when the slope of the repelling branch of the critical manifold is greater than the slope of the attracting branches of the critical manifold. As this is not the case in the Van der Pol system, it can be expected that only headless saddle-node canard cycles are possible in the Van der Pol system.

We finally point out that some quantitative information obtained in the manuscript could be relevant for applications. As example of such information we highlight the period of the canard cycles, see Lemma \ref{lem:times}, and the location of the saddle-node canards in terms of the parameter, see propositions \ref{prop:R3z} and \ref{prop:R4z}. The dependence between the height of a canard cycle and the bifurcation parameter $a$ at which it appears could be approximated from the estimation $|\tilde{a}-\hat{a}|$ appearing in Theorem \ref{th:exist_Gamma}.

\section{Acknowledgments}
The three authors want to thank M.  Desroches for introducing us into the problem and for his continuous support. They also thanks M. Krupa and S.  Rodrigues for their fruitful conversations.

First author is supported by Ministerio de Ciancia, Innovación y Universidades, through the project PGC2018-096265-B-I00. Second author is supported by Ministerio de Ciencia, Innovación y Universidades through the project RTI2018-093521-B-C31. Third author is supported by Ministerio de Economía y Competitividad through the
project MTM2017-83568-P (AEI/ERDF, EU). 

\appendix
\section{About Poincaré half-maps and times of flight}

In this section we summary the main technical results on piecewise linear dynamics that we use along the manuscript. In particular we provide approximations for the Poincaré half-map and for the time of flight of the solutions inside the regions of linearity.

Consider 
the Poincaré half-maps 
$\Pi_R$, $\Pi_L$ and $\Pi_{LL}$
having the switching lines $\{x=\sqrt{\varepsilon}\},$ $\{x=-\sqrt{\varepsilon}\}$ and $\{x=-1\}$, respectively, as crossing section. By using the contact points $\mathbf{p}_R,\mathbf{p}_L$ and $\mathbf{p}_{LL}$ and the vector field on these points $\dot{\mathbf{p}}_R,\dot{\mathbf{p}}_L$ and $\dot{\mathbf{p}}_{LL}$, respectively, the Poincaré half-maps 
can be parametrized by 
$\Pi_R(\mathbf{p}_R-u\dot{\mathbf{p}}_R)= \mathbf{p}_R+v\dot{\mathbf{p}}_R$, $\Pi_L(\mathbf{p}_L-u\dot{\mathbf{p}}_L)= \mathbf{p}_L+v\dot{\mathbf{p}}_L$, and  $\Pi_{LL}(\mathbf{p}_{LL}-u\dot{\mathbf{p}}_{LL})= \mathbf{p}_{LL}+v\dot{\mathbf{p}}_{LL}$, where $u$ and $v$ are positive values. Moreover, let ${\Pi}_{Ld}$ be the point transformation through the flow of the linear system $\dot{\mathbf{x}}=A_L\mathbf{x}+\mathbf{b}_L$, such that ${\Pi}_{Ld}(\mathbf{p}_{L}-u\dot{\mathbf{p}}_{L})= \mathbf{p}_{LL}-v\dot{\mathbf{p}}_{LL}$.

\begin{lemma}\label{lem:Pmaps}
Fix $\varepsilon_0>0$ small enough and $\nu \in (0,1)$, for $\varepsilon \in (0,\varepsilon_0)$ it follows that
\begin{align*}
\Pi_L^{-1}\left(\begin{pmatrix}-\sqrt{\varepsilon} \\ h\end{pmatrix}\right)- \mathbf{q}_0^L &\approx 
        \begin{pmatrix} 0 \\ h e^{-\frac {k h}{\varepsilon(\sqrt{\varepsilon}-a)}} \end{pmatrix}, \quad\text{for } h>(-m+\varepsilon^{-\nu}\lambda_L^s)(\sqrt{\varepsilon}+a),\\
\Pi_R\left(\begin{pmatrix}\sqrt{\varepsilon} \\ h\end{pmatrix}\right)-\mathbf{q}_1^R&\approx
        \begin{pmatrix} 0 \\ h e^{-\frac {h}{\varepsilon(\sqrt{\varepsilon}-a)}} \end{pmatrix},\quad\text{for } h>(m-\varepsilon^{-\nu}\lambda_R^s)(\sqrt{\varepsilon}-a),\\
\mathbf{q}_1^{LL}-\Pi_{LL}\left(\begin{pmatrix}-1 \\ h\end{pmatrix}\right) &\approx 
    \begin{pmatrix} 0 \\ (k-h) e^{-\frac {k(k-h)}{\varepsilon(1+a)}} \end{pmatrix},
    \quad\text{for } h< p_2+ \varepsilon^{-\nu} \lambda_{LL}^s(1+a),    \\
\mathbf{q}_0^L-{\Pi}_{Ld}^{-1}\left(\begin{pmatrix}-1 \\ h\end{pmatrix}\right) &\approx
        \begin{pmatrix} 0 \\ (k-h) e^{-\frac{k^2}{\varepsilon} \ln\left(\frac {1+a}{a+\sqrt{\varepsilon}}\right) } \end{pmatrix}, \quad\text{for }
        h\in(p_2-\lambda_L^q(1+a),p_2),\\
\end{align*}
where $p_2=k(1-\sqrt{\varepsilon})-m(\sqrt{\varepsilon}+a)$ is the second coordinate of the point $\mathbf{p}_{LL}$.
\end{lemma}
\begin{proof}
Next, we compute the expressions of the Poincaré half-maps $\Pi_R$ and ${\Pi}_{Ld}$, the remainder expressions in the lemma follows in a similar way. 

Following Chapter 3 in \cite{LT14}, the coordinates $u$ and $v$ of the points $\mathbf{p}$ and $\Pi_R(\mathbf{p})$ in the Krylov base $\{\mathbf{p}_R,\dot{\mathbf{p}}_R\}$, that is $\mathbf{p}=\mathbf{p}_R-u \dot{\mathbf{p}}_R$ and  $\Pi_{R}(\mathbf{p})=\mathbf{p}_R+v\dot{\mathbf{p}}_R$, are invariant through translations and linear transformation. Therefore, $u$ and $v$ satisfy that $v\in\left(0,-\frac {1}{\lambda_R^q}\right)$ and 
\[
\left(
 \frac {1+v\lambda_R^s}{1-u\lambda_R^s}
\right)^{\frac {\lambda_R^q}{\lambda_R^s}} =  \frac {1+v\lambda_R^q}{1-u\lambda_R^q},
\]
where $\lambda_R^s$ and $\lambda_R^q$ are the slow and the fast eigenvalues of the matrix $A_R$. For $\varepsilon$ small enough, $\lambda_R^s$ tends to zero and hence, taking the limit of the left side term in the previous identity as $\lambda_R^s$ tends to zero, we obtain that $v$ can be implicitly approximated by
\[
v=\frac 1{\lambda_R^q} \left(-1 + (1-u\lambda_R^q)e^{\lambda_R^q(u+v)}\right).
\]

Since the $u$-coordinate of $\mathbf{p}=(\sqrt{\varepsilon},h)^T$ in the Krylov base is 
\[
u=\frac {h-m(\sqrt{\varepsilon}-a)}{\varepsilon(\sqrt{\varepsilon}-a)}>0,
\]
we conclude that  
\begin{align*}
\Pi_R(\mathbf{p})
&=\mathbf{p}_R+\frac{1}{\lambda_R^q}\left(-1+(1-u\lambda_R^q)e^{\lambda_R^q(u+v)}\right)\dot{\mathbf{p}}_R
\\
&=\begin{pmatrix} 
    \sqrt{\varepsilon} \\ \\
    (m+\lambda_R^s)(\sqrt{\varepsilon}-a)-\left(\lambda_R^s(\sqrt{\varepsilon}-a) -u\varepsilon(\sqrt{\varepsilon-a})\right)e^{\lambda_R^q(u+v)}
\end{pmatrix}\\
&=\begin{pmatrix} 
    \sqrt{\varepsilon} \\ \\
    (m+\lambda_R^s)(\sqrt{\varepsilon}-a)-\left((m+\lambda_R^s)(\sqrt{\varepsilon}-a) -h\right)e^{\lambda_R^q(v+\frac {h-m(\sqrt{\varepsilon}-a)}{\varepsilon(\sqrt{\varepsilon}-a)})}
\end{pmatrix}\\
&=    \mathbf{q}_1^R+\begin{pmatrix}
    0 \\ \\
    \left(h-(m+\lambda_R^s)(\sqrt{\varepsilon}-a) \right)e^{\lambda_R^q v+\frac {h-m(\sqrt{\varepsilon}-a)}{\lambda_R^s (\sqrt{\varepsilon}-a)}}
\end{pmatrix}.
\end{align*}
Taking into account that $-1<\lambda_R^q v<0$, the exponent in the previous expression can be approximated by 
$\frac {h}{\lambda_R^s(\sqrt{\varepsilon}-a)}$ provided that
\[
\frac {h-m(\sqrt{\varepsilon}-a)}
{\lambda_R^s(\sqrt{\varepsilon}-a)}<-\varepsilon^{-\nu},
\]
for a $\nu$ value in $(0,1)$. From this we conclude that $h>(m-\varepsilon^{-\nu}\lambda_R^s)(\sqrt{\varepsilon}-a)$, what finish the proof of the the lemma. 

Consider now the transformation $\Pi_{Ld}$, from points of the form $\mathbf{p}_L-u\dot{\mathbf{p}}_L$ to points of the form $\mathbf{p}_{LL}-v\dot{\mathbf{p}}_{LL}$. Even when the relationship between the coordinates $u$ and $v$ of these points is not explicitly computed in \cite{LT14}, we can use the same arguments than there to obtain that $v\in \left(0,\frac {1}{\lambda_L^s}\right)$ and 
\[
\left( \frac {1-v\lambda_L^s}{1-u\lambda_L^s} \right)^{\frac {\lambda_L^q}{\lambda_L^s}}=
 r^{1-\frac {\lambda_L^q}{\lambda_L^s}} 
 \left( \frac {1-v\lambda_L^q}{1-u\lambda_L^q} \right),
\]
where $r=\frac{1+a}{\sqrt{\varepsilon}+a}$ satisfies that $\mathbf{p}_{LL}-\mathbf{e}_L=r(\mathbf{p}_{L}-\mathbf{e}_L)$. Taking the limit of the left side term in the previous identity as $\lambda_L^s$ tends to zero, we obtain a new implicit relation between the coordinates $u$ and $v$ given by
\[
u=\frac 1{\lambda_L^q} \left(1-r(1-v\lambda_L^q)e^{-\frac {\lambda_L^q}{\lambda_L^s}\ln(r)+\lambda_L^q (v-u)}\right).
\]
Since the $v$ coordinate of a point $(-1,h)=\mathbf{p}_{LL}-v\dot{\mathbf{p}}_{LL}$ satisfies that \[
v=\frac {k(1-\sqrt{\varepsilon})-m(\sqrt{\varepsilon}+a)-h}{\varepsilon(1+a)}\in \left(0,\frac 1{\lambda_L^s}\right),
\]
it follows that
\[
k(1-\sqrt{\varepsilon})-m(\sqrt{\varepsilon}+a)-\lambda_L^q(1+a)<h<k(1-\sqrt{\varepsilon})-m(\sqrt{\varepsilon}+a),
\]
and therefore, the preimage by ${\Pi}_{Ld}$ of the point $(-1,h)$ is 
\[
\mathbf{p}_L-u\dot{\mathbf{p}}_L=\mathbf{q}_0^L +
\begin{pmatrix}
0\\
\left(
(1+a)\lambda_L^s+m(a+\sqrt{\varepsilon})-k(1-\sqrt{\varepsilon})+h
\right)
e^{-1-\frac{\lambda_L^q(1+a) \ln(r)-k(1-\sqrt{\varepsilon})+m(\sqrt{\varepsilon}+a)+h}
{\lambda_L^s(1+a)}}
\end{pmatrix}.
\]
For $\varepsilon$ small  enough previous expression can be rewriten as 
\[
\mathbf{p}_L-u\dot{\mathbf{p}}_L=\mathbf{q}_0^L +
\begin{pmatrix}
0\\
\left(
(m+k)(a+\sqrt{\varepsilon})-\lambda_L^q(1+a)+h
\right)
e^{-\frac{k}{\varepsilon}
\ln\left(\frac{1+a}{\sqrt{\varepsilon}+a}\right)}
\end{pmatrix}.
\]
The lemma for $\Pi_{Ld}$ follows by considering $\varepsilon$ small enough.
\end{proof}

\begin{remark}\label{rem:interv}
From Lemma \ref{lem:Pmaps}, the preimage by $\Pi_L$ of a point is exponentially close to $\mathbf{q}_0^L$ if there exists $\nu\in(0,1)$ such that, its second coordinate, $h$, is at a distance  $d\geq\varepsilon^{-\nu}\lambda_L^s(\sqrt{\varepsilon}+a)$ from the point $\mathbf{p}_L$, in particular when $d> \lambda_L^s(\sqrt{\varepsilon}+a)$.
\end{remark}

Let $\tau_R(h)$ be the time of flight of the solution between the points $(\sqrt{\varepsilon},h)^T$ and $\Pi_R(\sqrt{\varepsilon},h)$, let $\tau_L(h)$ be the time of flight between $(-\sqrt{\varepsilon},h)^T$ and $\Pi_L(-\sqrt{\varepsilon},h)$, let $\tau_{LL}(h)$ be the time of flight between $(-1,h)^T$ and $\Pi_{LL}(-1,h)$, let ${\tau}_{Ld}(h)$ be the time of flight between $(-\sqrt{\varepsilon},h)^T$ and $\Pi_{Ld}(-\sqrt{\varepsilon},h)$ and let $\tau_{RR}$ be the time of flight between $\mathbf{q}_0^{RR}$ and $\mathbf{q}_1^R$, see Figure \ref{dib:geommetry}, that is, $\tau_{RR}=\tau_R(h_0)$ where $h_0$ is the second coordinate of $\mathbf{q}_1^{LL}$,

\begin{lemma} \label{lem:times}
For $\varepsilon>0$ and small enough it follows that
\begin{align*}
 \tau_R(h)&\approx-\frac 1{\lambda_R^s} 
  \ln\left(
      1 +  
      \frac {(m+\lambda_R^s)(\sqrt{\varepsilon}-a)-h}
      {(\lambda_R^q-\lambda_R^s)(\sqrt{\varepsilon}-a)}  
    \right),\\
 \tau_L(h)&\approx\frac 1{\lambda_L^s} 
  \ln\left(
      1 + 
      \frac {h+(m+\lambda_L^s)(\sqrt{\varepsilon}+a)}
      {(\lambda_L^q-\lambda_L^s)(\sqrt{\varepsilon}+a)}
    \right),\\
 {\tau}_{Ld}(h)&\approx\frac 1{\lambda_L^s} 
  \ln\left(
      1 + 
      \frac {h+m(\sqrt{\varepsilon}+a)+\lambda_L^s(2\sqrt{\varepsilon}+a-1)}
      {(\lambda_L^q-\lambda_L^s)(\sqrt{\varepsilon}+a)}  
    \right),\\
 \tau_{LL}(h)&\approx-\frac1{\lambda_{LL}^s} 
  \ln\left(
      1+ 
      \frac{ h+m(\sqrt{\varepsilon}+a)+k(\sqrt{\varepsilon}-1)+\lambda_{LL}^s(1+a)}
           { (\lambda_{LL}^q-\lambda_{LL}^s)(1+a) }
    \right),\\
    \tau_{RR}&\approx - \frac1{\lambda_R^s} 
  \ln\left(
      1+ 
      \frac{ \lambda_R^s(\sqrt{\varepsilon}-a)+\lambda_{LL}^s(1+a)+k(\sqrt{\varepsilon}-1)+2m\sqrt{\varepsilon} }
      { (\lambda_R^q-\lambda_R^s)(\sqrt{\varepsilon}-a)}
    \right).
\end{align*}
\end{lemma}
\begin{proof}
Consider a point $\mathbf{p}=(\sqrt{\varepsilon},h)^T$ and its image by the Poincaré map $\Pi_R(\mathbf{p})$. From Lemma \ref{lem:Pmaps}, since $\varepsilon$ is small enough, we can substitute the point $\Pi_R(\mathbf{p})$ by the exponentially close point $\mathbf{q}_1^R$. In an equivalent way, we approximately compute the time $\tau_R(h)$ as the time of flight of the solution for travelling from $\mathbf{p}$ to $\mathbf{q}_1^R$. To do that we project the point $\mathbf{p}$ onto the point $\mathbf{p}_s$ contained in the slow manifold $\mu_R$, see \eqref{def:Se}, by following the fast eigenvector, $\mathbf{v}_R^q$. The point $\mathbf{p}_s$ is obtained by solving with respect the unknowns $r_s,r_q\in \mathbb{R}^+$ the linear system of equations 
\[
\mathbf{p}_s=\mathbf{p}-r_q\mathbf{v}_R^q=\mathbf{e}_R-r_s \mathbf{v}_R^s.
\]
We conclude that 
\[
r_s=\frac{(\sqrt{\varepsilon}-a)(m-1-\lambda_R^s)-h}{\lambda_R^s(\lambda_R^s-\lambda_R^q)} 
=\frac{(\sqrt{\varepsilon}-a)(m+\lambda_R^q)-h}{\lambda_R^s(\lambda_R^s-\lambda_R^q)}.
\]

Then, we compute $\tau_R(h)$ as the time of flight of the solution to travel from the projected point $\mathbf{p}_s$ to $\mathbf{q}_1^R$, that is
\[
e^{\lambda_R^s \tau_R(h)}=\frac {\|\mathbf{q}_1^R-\mathbf{e}_R\|}{\|\mathbf{p}_s-\mathbf{e}_R\|}
=\frac{\frac {a-\sqrt{\varepsilon}}{\lambda_R^s} \|\mathbf{v}_R^s\|}{r_s \|\mathbf{v}_R^s\|}
=\frac {(\sqrt{\varepsilon}-a)(\lambda_R^q-\lambda_R^s)}{(\sqrt{\varepsilon}-a)(m+\lambda_R^q)-h}.
\]
The lemma follows by isolating $\tau_R(h)$. The remainder of the functions are computed by following similar arguments.
\end{proof}

\begin{remark}
In the computation of function $\tau_R(h)$ there is only one approximation, which consists in using the exponentialy close point $\mathbf{q}_1^R$, instead of the point $\Pi_R(\mathbf{p})$.
\end{remark}

For the purposes of this manuscript it is also important the knowledge of the map $\Phi(x_0)$ defined by the flow through the points on the fast nullcline $(x_0,f(x_0))$. Since the nullcline is piecewise, so is the map $\Phi$, which domain depends on $x_0<-1$ and $x_0>-1$. In particular, when $x_0\in (x_r,x_u)$ the map $\Phi$ is defined as $\Phi_{4z}(x_0)$ by the second coordinate of the point on the switching line $\{x=-1\}$ defined backward through the flow, and when $x_0\in [-1,x_s)$ the map $\Pi$ is defined as $\Phi_{3z}(x_0)$ by the second coordinate of the point on the switching line $\{x=-\sqrt{\varepsilon}\}$ defined forward through the flow. In the next result we relate the map $\Phi$ with the Poincaré maps analysed in Lemma \ref{lem:Pmaps}. This relationship is established through the values $x_r$, $x_s$ and $x_u$ introduced in \eqref{def:xs} and \eqref{def:xu}, respectively.

\begin{lemma}\label{lem:dom}
For $\varepsilon>0$ and $\Phi(x_0)=h$ it follows that:
\begin{itemize}
    \item [a)] If $x_0\in(x_r,x_u)$ there exists $\nu\in (0,1)$ such that $h_r<h<h_u$, where 
    $h_r=-(m+\lambda_L^s)(\sqrt{\varepsilon}+a)$ and $h_u=-m(\sqrt{\varepsilon}+a)+k(1-\sqrt{\varepsilon}).
    +\varepsilon^{-\nu}\lambda_{LL}^s(1+a),$
    \item [b)] If $x_0\in[-1,x_s)$ there exists $\nu\in (0,1)$ such that $h_s<h<h_M$, where $h_s=-(m-\varepsilon^{-\nu}\lambda_{L}^s)(\sqrt{\varepsilon}+a)$ and $h_M=-m(\sqrt{\varepsilon}+a)+k(1-\sqrt{\varepsilon}).$
\end{itemize}
\end{lemma}
\begin{proof}
We restrict ourselves to the proof of statement (b). The statement (a) of the lemma follows by using similar arguments.

Let $r$ and $v$ be the coordinates of the points $(x_0,f(x_0))$ and $(-\sqrt{\varepsilon},\Phi(x_0))$, respectively, in the  Krylov base $\{\mathbf{p},\dot{\mathbf{p}}\}$ with $\mathbf{p}=\mathbf{p}_L-\mathbf{e}_L$. Thus, there exists $\tau_f>0$  such that $e^{A_L \tau_f}r\mathbf{p}=\mathbf{p}+v\dot{\mathbf{p}}$,where 
\[
r=\frac {a-x_0}{a+\sqrt{\varepsilon}},
\]
is obtained from $(x_0,f(x_0))^T-\mathbf{e}_L=r(\mathbf{p}_L-\mathbf{e}_L)$.

Following similar arguments than those in Chapter 3 in \cite{LT14}, we obtain next relation between the coordinates $r$ and $v$,
\begin{equation}\label{aux:lem_dom}
{r^{\lambda_L^s-\lambda_L^q}}(1+v\lambda_L^s)^{\lambda_L^q}=
(1+v\lambda_L^q)^{\lambda_L^s},
\end{equation}
 In particular, when $x_0=x_s$, it follows that $r=1+\lambda_L^s$ and the corresponding coordinate $v$ satisfies
\begin{equation}\label{aux:lem3}
e^{\lambda_L^s-\lambda_L^q}=e^{-v\lambda_L^q}(1+v\lambda_L^q),
\end{equation}
where we have taken into account that $(1+v\lambda_L^s)^{\frac {\lambda_L^q}{\lambda_L^s}}$ and $(1+\lambda_L^s)^{\frac {\lambda_L^s-\lambda_L^q}{\lambda_L^s}}$ tend to $e^{v\lambda_L^q}$ and $e^{\lambda_L^s-\lambda_L^q}$, respectively, as $\lambda_L^s$ tends to zero.

From expression \eqref{aux:lem3} and for $\varepsilon$ small enough it follows that $v>1$ and, hence, $\mathbf{p}+v\dot{\mathbf{p}}=(-\sqrt{\varepsilon},h)^T-\mathbf{e}_L$ with $h>-m(\sqrt{\varepsilon}+a)+\varepsilon(\sqrt{\varepsilon}+a)$. Therefore, for a fixed $\varepsilon>0$ there exists $\nu\in(0,1)$ such that $h>-m(\sqrt{\varepsilon}+a)+\varepsilon^{-\nu}\lambda_L^s(\sqrt{\varepsilon}+a)$, which proves the lower bound in the statement (b).

For the upper bound, we note that when $x_0=-1$ then $r=\frac {a+1}{a+\sqrt{\varepsilon}}$, which tends to infinity as $\varepsilon^{-1/2}$ when $\varepsilon$ tends to zero. From the relation  between the variables $v$ and $r$, see \eqref{aux:lem_dom}, we conclude that for $\varepsilon$ small enough, expression $\frac {1+v\lambda_L^s}{r}$ tends to 1 if $\varepsilon$ tends to 0. Therefore $v$ tends  to 
\[
v=\frac {1-r}{\lambda_L^s}=\frac 1{\lambda_L^s}\frac {1-\sqrt{\varepsilon}}{a+\sqrt{\varepsilon}}.
\]
Hence, since $\mathbf{p}+v\dot{\mathbf{p}}=
(\sqrt{\varepsilon},h)$, it follows that 
\[
h=-m(\sqrt{\varepsilon}+a)+v \varepsilon(a+\sqrt{\varepsilon})\approx-m(\sqrt{\varepsilon}+a)+k(1-\sqrt{\varepsilon}). 
\]
This ends the proof of the statement (b).
\end{proof}

\end{document}